\documentclass{article}

\usepackage{geometry}
\usepackage{graphicx}
\usepackage{multicol}
\usepackage[bottom]{footmisc}
\usepackage[binary-units=true]{siunitx}

\usepackage{pgfplots}
\usepackage{pgfplotstable}
\usepackage{booktabs}
\usepackage{authblk}
\usepackage{amsmath,amsthm,mathtools}
\usepackage{amsfonts,amssymb}
\usepackage{hyperref}
\usepackage{float} %
\usepackage[ruled, longend, linesnumbered, noline]{algorithm2e}
\SetEndCharOfAlgoLine{}

\newcommand{\restrictedto}[2]{{%
		\left.\kern-\nulldelimiterspace %
		{#1} %
		\vphantom{\big|} %
		\right|_{#2} %
}}

\newcommand{\Cheb}{\operatorname{Cheb}}
\newcommand{\diag}{\operatorname{diag}}
\renewcommand{\div}{\operatorname{div}}

\newcommand{\vr}[1]{\mathbf{#1}}

\newtheorem{Lemma}{Lemma} 
\newtheorem{Theorem}{Theorem} 

\begin{document}

\title{Matrix-free Monolithic Multigrid Methods\\ for Stokes and Generalized Stokes Problems
\footnote{This work has been supported by the Austrian Science Fund (FWF) grant P29181 `Goal-Oriented Error Control for Phase-Field Fracture Coupled to Multiphysics Problems', 
and by the Doctoral Program W1214-03 at the Johannes Kepler University Linz.
The third author has been partially funded by the Deutsche Forschungsgemeinschaft (DFG) under Germany’s Excellence Strategy within the Cluster of Excellence PhoenixD (EXC 2122, Project ID 390833453).
}}

\renewcommand\Affilfont{\itshape\small}

\author[1]{{D. Jodlbauer}}
\author[1,2]{{U. Langer}}
\author[3,4]{{T. Wick}}
\author[1]{{W. Zulehner}}

\affil[1]{Johann Radon Institute for Computational Mathematics, Austrian Academy of Sciences, Altenbergerstr. 69, A-4040 Linz, Austria}
\affil[2]{Institute for Computational Mathematics, Johannes Kepler University Linz, Altenbergerstr. 69, A-4040 Linz, Austria}
\affil[3]{Leibniz Universit\"at Hannover, Institut f\"ur Angewandte Mathematik, Welfengarten 1, 30167 Hannover, Germany}
\affil[4]{Cluster of Excellence PhoenixD (Photonics, Optics, and
	Engineering - Innovation Across Disciplines), Leibniz Universit\"at Hannover, Germany}

\maketitle

\abstract{
We consider the widely used continuous 
$\mathcal{Q}_{k}$-$\mathcal{Q}_{k-1}$
quadrilateral or hexahedral  Taylor-Hood elements
for the finite element discretization of the 
Stokes and generalized Stokes systems in two and three spatial dimensions.
For the fast solution of the corresponding symmetric, but indefinite system of finite element equations, 
we propose and analyze matrix-free monolithic geometric multigrid 
solvers 
that are based on appropriately scaled Chebyshev-Jacobi smoothers.  
The analysis is based on results 
by Sch\"oberl and Zulehner (2003). We present and discuss  several numerical results for typical
benchmark problems.\\
}
\section{Introduction}
\label{sec:Introduction}

The Stokes problem 
\begin{equation}
\label{eqn:Stokes:cf}
 - \mu \Delta u + \nabla p =  f, \; \div u = 0 \; \mbox{in}\, \Omega,\;
 \mbox{and}\, u=u_D:=0 \; \mbox{on}\, \partial \Omega
\end{equation}
and the generalized Stokes problem 
\begin{equation}
\label{eqn:genStokes:cf}
 - \mu \Delta u 
+ \gamma \varrho u + \nabla p =  f, \; \div u = 0 \; \mbox{in}\, \Omega,\;
 \mbox{and}\, u=u_D:=0 \; \mbox{on}\, \partial \Omega
\end{equation}
are two important basic problems in fluid mechanics, where one looks for 
the velocity field $u=(u_1,\dots,u_d)$ and pressure field $p$ for some 
given right-hand side $f=(f_1,\dots,f_d)$, dynamic viscosity $\mu$, and
density $\varrho$.
For simplicity, we consider homogeneous Dirichlet date $u_D$ for the velocity 
in the analysis part of the paper, whereas non-homogeneous Dirichlet date $u_D$ 
are also permitted in our numerical studies presented in Section~\ref{sec:NumericalResults}.
Here, $d \in \{2,3\}$ denotes the space dimension.
The non-negative parameter $\gamma$ depends on 
the time integration scheme, e.g., implicit Euler gives $\gamma =1 / \Delta t$,
and $\gamma = 0$ formally turns the generalized Stokes problem \eqref{eqn:genStokes:cf}
into the classical Stokes problem \eqref{eqn:Stokes:cf}.
The computational domain $\Omega \subset \mathbb{R}^d$ is supposed to be 
bounded and Lipschitz with boundary $\partial \Omega$.
We mention that
the Stokes problem arise from the iterative linearization
of the stationary Navier-Stokes equations when treating the 
advection term
on the right-hand side. The same idea leads to the generalized Stokes problem 
after an implicit time discretization of the instationary Navier-Stokes equations;
see, e.g., \cite{Turek:1999a}.
Moreover, both problems are building blocks in fluid-structure interaction (FSI) 
problems where the Navier-Stokes equations model the fluid part;
see, e.g., \cite{Richter:2017a}, \cite{Wichrowski:2021a}, \cite{Jodlbauer:2021a}.
The Stokes problem itself can be used as model for describing dominantly viscous fluids. 
Thus fast solvers for the discrete surrogates of \eqref{eqn:Stokes:cf} or \eqref{eqn:genStokes:cf}
are of great importance in many applications 
where the Navier-Stokes equations are involved.
The mixed finite element discretization of \eqref{eqn:Stokes:cf} or \eqref{eqn:genStokes:cf}
finally leads to large-scale linear systems of the form 
\begin{equation}
\label{eqn:LinearSystem}
\begin{pmatrix} A & B^T\\ 
                B & -C 
\end{pmatrix}
\begin{pmatrix} \vr u \\ \vr p \end{pmatrix} = \begin{pmatrix} \vr f\\\vr g \end{pmatrix},
\end{equation}
where here $\vr u \in \mathbb{R}^n$, $\vr p\in \mathbb{R}^m$, $\vr f\in \mathbb{R}^n$, and $\vr g\in \mathbb{R}^m$ are vectors, 
$A$ is a symmetric and positive definite (spd) $n \times n$ matrix, 
$C$ is a non-negative $m \times m$ matrix, $B^T$ is transposed to the $m \times n$  matrix $B$. 
Usually the matrix $C$ is $0$, but stabilizations can create a non-trivial matrix $C$;
see, e.g., \cite{ElmanSilvesterWathen:2005a,John:2016a}.
 
The system matrix of \eqref{eqn:LinearSystem} is symmetric, but indefinite.
So special direct or iterative solvers tailored to these properties of the system matrix 
are needed; see, e.g., \cite{Saad:2003a}.
Iterative methods are most suited for large-scale systems as they usually arise after 
mixed finite element discretization. 
Similar to the Conjugate Gradient (CG) method in the case of spd systems, 
the MinRes \cite{PaigeSaunders:1975a}
and the Bramble-Pasciak CG \cite{BramblePasciak:1988a}
are certainly 
efficient Krylov-subspace solvers for symmetric, indefinite systems of the form \eqref{eqn:LinearSystem}
provided that suitable preconditioners are available. 
We refer the reader to the survey article \cite{BenziGolubLiesen:2005a},
the monograph \cite{ElmanSilvesterWathen:2005a}, 
the more recent articles \cite{Notay_2014,Axelsson_2015,Notay_2019}, 
and the references therein.
\\
Monolithic geometric multigrid methods are another class of efficient solvers
for large scale systems of finite element equations.
The multigrid convergence analysis is based on the {\it smoothing property} 
and the {\it approximation property}; see \cite{Hackbusch:1985a}.
In the case of symmetric, indefinite systems of the form \eqref{eqn:LinearSystem},
the construction of suitable smoothers is crucial for the overall
efficiency of the multigrid method.
{
Such smoothers range  from %
low cost methods like Richardson-type iterations %
applied to the normal equation associated with \eqref{eqn:LinearSystem}, see \cite{Verfuerth1984,
Brenner_1996}, 
collective smoothers, like Vanka-type smoothers, based on the solutions of small local problems, see \cite{Vanka_1986,
SchoeberlZulehner:2003a}, 
to more expensive methods, 
whose costs for the Stokes problem are comparable with the costs of an optimal preconditioner for a discrete Laplace-like equation either in $\mathbf{p}$, like the exact and inexact Braess-Sarazin smoothers
\cite{BraessSarazin:1997a,Zulehner:2000a} or, more recently in $\mathbf{u}$, resulting in better smoothing properties \cite{Brenner_2014, Brenner_2017}. 
A general framework for the construction and the analysis of smoothers was introduced in \cite{Wittum:1989a,Wittum:1990a}
with the concept of transforming smoothers, which include several classes like distributive Gauss-Seidel smoothers.
Yet another class of smoothers 
are Uzawa-type smoothers, see the survey article \cite{DrzisgaJohnRuedeWohlmuthZulehner:2018a} (as well as the references therein), where also other classes of block smoothers are discussed.  %
}

In this paper, we are looking for an efficient (parallel) matrix-free implementation 
of all ingredients of  Geometric MultiGrid (GMG) algorithms, 
which are based on a solid convergence analysis.
Current developments of matrix-free solvers include 
problems in finite‐strain hyperelasticity \cite{DaPeAr20}, 
phase-field fracture  \cite{JoLaWi20_CMAME,JoLaWi20,Jodlbauer:2021a},
fluid-structure interaction \cite{Wichrowski:2021a},
discontinuous Galerkin \cite{KrKo19},
compressible Navier-Stokes equations \cite{guermond2021implementation}, incompressible 
Navier-Stokes and Stokes equations \cite{FrCaAnPa20} as well as 
sustainable open-source code developments \cite{MuKoKr20,CleHeiKaKr21}, 
matrix-free implementations on locally refined meshes \cite{MuHeiSaaKr22_arxiv},
and implementations on graphics processors \cite{KroLjun19}.

Matrix-free implementations should avoid the assembling of all matrices involved  
in the multigrid algorithm. Matrix-free implementations can directly make use 
of the finite-difference-star representation of a matrix-vector multiplication 
especially, for finite difference or volume discretizations, or similarly
the stencil representation of finite elemenet equations \cite{BaDrMo18},
or can be realized by 
the element matrices on the fly in the case of finite element discretizations; see, e.g., 
\cite{KronbichlerKormann:2012a}.
We are going to use the latter (element-based) matrix-free technique that leads not only to a considerable 
reduction of the memory demand, but also to a reduction of the arithmetical 
complexity, especially, for higher-order finite element discretizations;
see, e.g., \cite{MuKoKr20}.
Moreover, 
element-based matrix-free techniques are very suited for the implementation
on massively parallel computers with distributed memory \cite{KronbichlerKormann:2012a}.
On the other side, the element-based matrix-free technology 
restricts the use of smoothers mentioned above. 
Smoothers that are constructed from the assembled matrix are out of the game.
We will consider inexact symmetric Uzawa smoothers {(also known as a class of block approximate factorization smoothers)}  of the form
\begin{align}
\label{eqn:isus1}
	\hat{\vr u}^{j+1}     &= \vr u^j + {\hat A}^{-1} (\vr f - A\vr u^j -B^T \vr p^j), \\
\label{eqn:isus2}	
	\vr p^{j+1}           &= \vr p^j + {\hat S}^{-1} (B \hat{\vr u}^{j+1} - C\vr p^j - \vr g), \\
\label{eqn:isus3}	
	\vr u^{j+1}           &= \vr u^j + {\hat A}^{-1} (\vr f - A\hat{\vr u}^{j+1} - B^T \vr p^j),
\end{align}
which were analyzed by Sch\"oberl and Zulehner in \cite{SchoeberlZulehner:2003a}.
They showed the smoothing property provided that the spd ``preconditioning'' matrices $\hat A$ and $\hat S$
satisfy the spectral inequalities
\begin{equation}
\label{eqn:Smoother:SpectralInequalities}
  {\hat A} \ge A \quad \mbox{and}\quad  {\hat S} \ge {\tilde S} = B{\hat A}^{-1} B^T + C,
\end{equation}
and an additional estimate of the difference of the system matrix and the ``preconditioning'' 
matrix generated by the smoother \eqref{eqn:isus1} - \eqref{eqn:isus3}.
The W-cycle multigrid convergence then follows from this smoothing property 
and the approximation property.
The application of  ${\hat A}^{-1}$ and ${\hat S}^{-1}$ to vectors will be realized completely 
matrix-free by means of the Chebyshev-Jacobi method. 
Once the smoother can be performed matrix-free,
the complete multigrid cycle can be implemented matrix-free and in parallel.
The implementation is based on the deal.II finite element library \cite{dealii}, and in particular, the matrix-free and geometric multigrid modules 
\cite{KronbichlerKormann:2012a,JaKa11}.
The parallel performance and scalability of matrix-free GMG solvers has 
widely 
been 
demonstrated; see, e.g., \cite{KrKo19,MuKoKr20,DaPeAr20,JoLaWi20,MuHeiSaaKr22_arxiv}.
Hence, in this work, we mainly focus  on the numerical analysis and  
the quantitative numerical illustration of the theoretical results.

The remainder of the paper is organized as follows. 
Section~\ref{sec:mVF+mFEM} provides the mixed variational formulation of 
\eqref{eqn:Stokes:cf} and \eqref{eqn:genStokes:cf}, their mixed fined element 
discretization by means of inf-sup stable $\mathcal{Q}_{k}$-$\mathcal{Q}_{k-1}$ quadrilateral or hexahedral finite 
elements, and recalls some well-known results on solvability and 
discretization error estimates.  
In Section~\ref{sec:MatrixfreeGMG}, we present the matrix-free GMG algorithm, in particular,
the class of matrix-free smoothers we are going to use, and the analysis of the 
smoothing and approximation properties yielding W-cycle multigrid convergence.
Section~\ref{sec:NumericalResults} presents and discusses our numerical results for different
benchmark problems 
supporting
the theoretical results. %
Finally, we draw some conclusions and give an outlook in Section~\ref{sec:ConclusionsOuilook}.

\section{Mixed Variational Formulation and Discretization}
\label{sec:mVF+mFEM}
The mixed variational formulation of \eqref{eqn:Stokes:cf} or \eqref{eqn:genStokes:cf} read as follows: 
Find $u \in V = H^1_0(\Omega)^d$ and 
$p \in Q = L_0^2(\Omega) = \{q \in L^2(\Omega): \int_\Omega q dx = 0\}$ such that
\begin{alignat}{3}
\label{eqn:mvf1}
	a(u,v) + b(v,p)    &= \langle F,v \rangle &\quad& \forall v \in V, \\
\label{eqn:mvf2}	
	b(u,q)             &= 0                   && \forall q \in Q, 
\end{alignat}
with
\begin{align}
\label{eqn:a(u,v)}
	a(u,v) &= \int_\Omega (\mu\, \nabla u : \nabla v + \gamma \varrho u v) \,  dx,\\
\label{eqn:b(v,q)}	
	b(u,q) &=  \int_\Omega q \, \div u \, dx,\\
\label{eqn:Fv}	
\langle F,v \rangle  &= \int_\Omega f v\, dx.
\end{align}
In some applications, e.g., FSI, instead of the bilinear form \eqref{eqn:a(u,v)},
the bilinear form
\begin{equation}
\label{eqn:a(varepsilon(u),varepsilon(v))}
	a(u,v) = \int_\Omega (\mu\, \varepsilon(u) : \varepsilon(v) + \gamma \varrho u v) \,  dx,
\end{equation}
is used, where $\varepsilon(u) = \frac{1}{2} (\nabla u +  (\nabla u)^T)$ denotes 
the deformation rate tensor, and $\mu$ and $\varrho$ may depend on the spatial 
variable $x$, but should always be uniformly positive and bounded;
see, e.g., \cite{Wichrowski:2021a}.
In both cases, the bilinear form $a(\cdot,\cdot)$ is elliptic and bounded on $V$.
Since $F$ obviously belongs to the dual space $V^*$ of $V$, and 
since the bilinear form $b(\cdot,\cdot)$ is also bounded and fulfills
the famous LBB (inf-sup) condition
\begin{equation}
	\label{eqn:infsup}
	\inf_{q \in Q \setminus \{0\}} \sup_{v \in V \setminus \{0\}} \frac{b(v,q)}{\|v\|_{H^1(\Omega)^d} \|q\|_{L^2(\Omega)}}  \ge c  > 0,	
\end{equation}
the mixed variational problem \eqref{eqn:mvf1} - \eqref{eqn:mvf2} is well-posed,
i.e. there exists a unique solution $(u,p) \in V\times Q$ satisfying an 
a priori estimate; see, e.g., \cite{GiRa1986,Temam2001}.
Throughout this paper, $c$ denotes a generic constant independent of the mesh level.

For later use we mention that the mixed variational problem \eqref{eqn:mvf1} -  \eqref{eqn:mvf2} can also be written as a variational problem on $V \times Q$:
Find $(u,p) \in V \times Q$ such that
\begin{equation} \label{mixed}
  \mathcal{B}((u,p),(v,q)) =  \langle F,v \rangle 
  \quad \text{for all} \ (v,q) \in V \times Q
\end{equation}
with the bilinear form
\[
  \mathcal{B}((w,r),(v,q)) = a(w,v) + b(v,r) + b(w,q) .
\]

\begin{figure}[htb]
	\centering
	\includegraphics{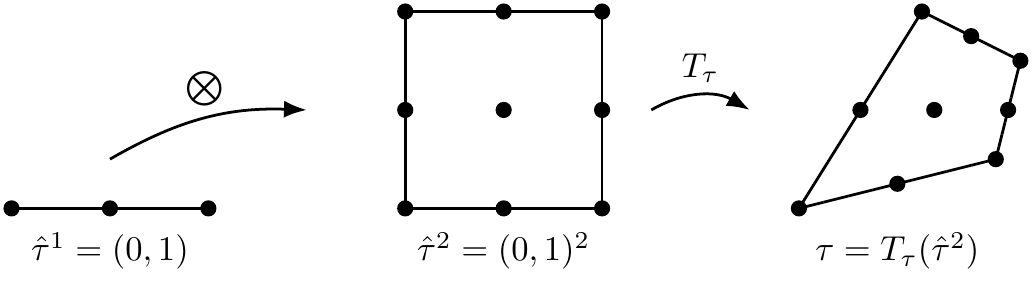}
	\caption{Construction of the tensor product nodes using a bilinear mapping $T_\tau$.}
	\label{fig:support points}
\end{figure}

\begin{figure}[htb]
	\centering
	\includegraphics[width=0.7\textwidth]{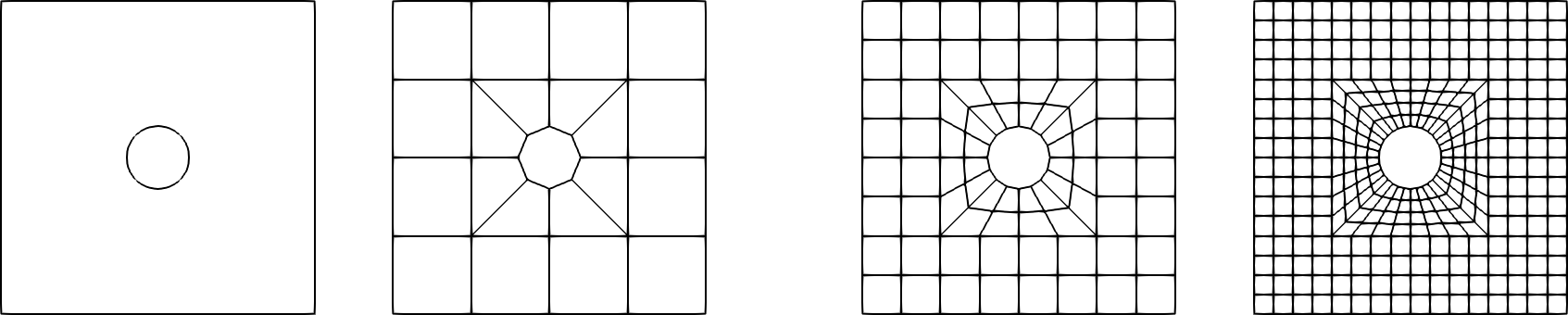}
	\caption{The mesh is gradually adapted to the geometry (left) during refinement.}
	\label{fig:mesh}
\end{figure}

Let us now briefly describe the mixed finite element discretization of \eqref{eqn:mvf1} - \eqref{eqn:mvf2}
by means of 
Taylor-Hood
$\mathcal{Q}_{k}$-$\mathcal{Q}_{k-1}$ finite elements, i.e., continuous, piecewise polynomials of degree $k$ and $k-1$ in each coordinate direction, with $k \in \mathbb{N}, k \geq 2$.
We assume that the computational domain $\Omega$ can be decomposed into 
non-overlapping 
convex quadrilateral (hexahedral)
elements $\tau \in \mathcal{T}_h$.
Each of the elements $\tau \in \mathcal{T}_h$ is the image of the reference element $(0,1)^d$ 
by a bilinear (trilinear) transformation $T_\tau$.
The nodal basis points for the shape functions are constructed using a tensor-product scheme on the reference element $(0,1)^d$ as visualized in Figure~\ref{fig:support points}.
This leads to the elementwise space
\[
	\mathcal{Q}_{k}(\tau) := \{ \hat{u} \circ T_\tau, \hat{u} \in \mathcal{Q}_{k} \},
\]
and the corresponding global, continuous space
\[
	V_h := \{ v \in V : \restrictedto{v}{\tau} \in \mathcal{Q}_{k}(\tau) \quad \forall \tau \in \mathcal{T}_h \}.
\]
The same construction is applied for the pressure space to obtain
\[
	\mathcal{Q}_h := \{ q \in Q \cap C(\Omega) : \restrictedto{q}{\tau} \in \mathcal{Q}_{k-1}(\tau) \quad \forall \tau \in \mathcal{T}_h \}.	
\]
Since not all features of the geometry can be represented by this construction (e.g., circles), we adapt the mesh during refinement to resolve those entities; see Figure~\ref{fig:mesh}.

Now the mixed finite element discretization of \eqref{eqn:mvf1} - \eqref{eqn:mvf2} 
reads as follows:
Find $u_h \in V_h \subset V$ and $p_h \in Q_h \subset Q $ such that
\begin{align}
\label{eqn:mfem1}
	a(u_h,v_h) + b(v_h,p_h)    &= \langle F,v_h \rangle \quad \forall v_h \in V_h, \\
\label{eqn:mfem2}	
	b(u_h,q_h)             &= 0 \quad \forall v_h \in Q_h. 
\end{align}
Once the usual nodal basis is chosen in $V_h$ and $Q_h$, the finite element scheme 
\eqref{eqn:mfem1} - \eqref{eqn:mfem2} is equivalent to the linear system \eqref{eqn:LinearSystem} 
with $C = 0$ and $g=0$ that can be written in the compact form
\begin{equation}
\label{eqn:LinearSystemKx=f}
\mathbf{K} \mathbf{x} = \mathbf{b},\quad \mbox{with}\; 
\mathbf{K} = \begin{pmatrix} A & B^T\\ B & 0 \end{pmatrix},\; 
\mathbf{x} = \begin{pmatrix} \vr u \\ \vr p \end{pmatrix},
\mbox{ and }
\mathbf{b} = \begin{pmatrix} \vr f\\0 \end{pmatrix}.
\end{equation}

The discrete inf-sup condition
\begin{equation}
\label{eqn:disinfsup}
\inf_{q_h \in Q_h \setminus \{0\}} \sup_{v_h \in V_h \setminus \{0\}} \frac{b(v_h,q_h)}{\|v_h\|_{H^1(\Omega)^d} \|q_h\|_{L^2(\Omega)}}  \ge c > 0,	
\end{equation}
ensures existence and uniqueness of the finite element solution $(u_h,p_h) \in V_h \times Q_h$,
and the estimate of the discretization error by the best approximation error. 
The inf-sup condition is ensured for $d=2$ and all $k \ge 2$, see \cite{Stenberg1990}. The inf-sup condition also holds
for $d=3$ and all $k \ge 2$ under mild conditions on the mesh, if we restrict ourselves to affine linear transformations $T_\tau$.
We refer the reader to Appendix A for the precise formulation of the conditions imposed on the mesh.
As a consequence the following Cea-type discretization error estimate holds.
\[
  \|u - u_h\|_{H^1(\Omega)^d} + \|p - p_h\|_{L^2(\Omega)} \le c \, 
  \left( \inf_{v_h \in V_h} \|u - v_h\|_{H^1(\Omega)^d} + \inf_{q_h \in Q_h} \|p - g_h\|_{L^2(\Omega)} \right).
\]
If $\Omega$ is convex or if its boundary is sufficiently smooth, it follows from the theory of elliptic regularity that $(u,p)$ even belongs to $H^2(\Omega)^d \times H^1(\Omega)$ and we have
\begin{equation}
  \|u\|_{H^2(\Omega)^d} + \|p\|_{H^1(\Omega)} \le c \, \|f\|_{L^2(\Omega)}.
\end{equation}
Based on this elliptic regularity and a standard duality argument we additionally have the following error estimate in $L^2(\Omega)^d$:
\[
  \|u - u_h\|_{L^2(\Omega)^d} \le c \, h \, \left( \|u - u_h\|_{H^1(\Omega)^d} + \|p - p_h\|_{L^2(\Omega)}\right),
\]
where $h$ denotes the mesh size of $\mathcal{T}_h$;
see, e.g., \cite{John:2016a}.

\section{Matrix-free Geometric Multigrid Solvers}
\label{sec:MatrixfreeGMG}

\subsection{Multigrid Algorithm and its Matrix-Free Implementation}
\label{subsec:MultigridAlgorithm+MatrixfreeImplementation}

Now let us assume that there is a sequence of meshes 
$\mathcal{T}_\ell$, $\ell=0,1,...,L$, where the subindex $0$
denotes the coarsest mesh,
and the subindex $L$ the finest mesh obtained after $L$ uniform 
mesh refinement steps.
This sequence of meshes leads to a sequence of nested finite element spaces $V_\ell \times Q_\ell$.

The multigrid method starts on the finest level $L$, and recursively adapts the current solution by solving the equation on the coarser levels.
The corrections are then prolongated back to the finer levels.
A crucial component is the smoother, which is applied before and after switching levels.
{Sometimes the post-smoothing is omitted.}
Its main purpose is to "smooth" the error, such that it can be represented on the coarser grids without losing too much accuracy.
The scheme is illustrated in Algorithm~\ref{alg:mg}.
More details and the analysis for multigrid methods can be found in the standard literature; 
see, e.g., \cite{Hackbusch:1985a,Br93}.

\begin{algorithm}
{
	\SetKwFunction{MG}{MG}
	\SetKwInOut{Input}{Input}	
	
	\MG($\mathbf{x}_\ell, \mathbf{b}_\ell$)
	
	\Indp

	\If{$\ell = 0$}
	{
		$\mathbf{x}_\ell \gets \mathbf{K}^{-1}_\ell \mathbf{b}_\ell$\tcc*{coarse grid solution}
	}
	
	$\mathbf{x}_\ell \gets \mathcal{S}_\ell(\mathbf{x}_\ell, \mathbf{b}_\ell)$ \tcc*{pre-smoothing}
	
	$\mathbf{r}_\ell \gets \mathbf{b}_\ell - \mathbf{K}_\ell \mathbf{x}_\ell$ \tcc*{residual}
	
	$\mathbf{r}_{\ell - 1} \gets \mathbf{I}^{\ell-1}_{\ell} \mathbf{r}_\ell $ \tcc*{restriction}\label{alg:mg:restriction}
	
	$\mathbf{w}_{\ell - 1} \gets 0$
	
	\For{$j=1, \dots, \gamma$}
	{
		$\mathbf{w}_{\ell - 1} \gets \MG(\mathbf{w}_{\ell - 1}, \mathbf{r}_{\ell - 1})$\tcc*{recursion}
	}
	
	$\mathbf{w}_\ell \gets \mathbf{I}_{\ell-1}^{\ell} \mathbf{w}_{\ell - 1}$\tcc*{prolongation}
	
	$\mathbf{x}_\ell \gets \mathbf{x}_\ell + \mathbf{w}_\ell$\tcc*{correction}
	
	$\mathbf{x}_\ell \gets \mathcal{S}_\ell(\mathbf{x}_\ell, \mathbf{b}_\ell)$\tcc*{post-smoothing}\label{alg:mg:post-smoothing}
	
	\Return{$\mathbf{x}_\ell$}
	
	\caption{Multigrid method for the solution of $\mathbf{K}_\ell \mathbf{x}_\ell = \mathbf{b}_\ell$.}
	\label{alg:mg}
}
\end{algorithm}

GMG methods are particularly suited for matrix-free applications
{since the entries of the matrix $\mathbf{K}=\mathbf{K}_\ell$ 
(we omit the level subindex where it is not needed)
are actually not required.}
{Therefore, we would like to avoid}
the assembly of the huge sparse matrix 
{$\mathbf{K}$} %
and its subsequent multiplications with a vector 
{$\mathbf{x}$}. %
Rather, we assemble the result
{$\mathbf{Kx}$} %
directly, i.e.,
\[
{
\mathbf{Kx} = \sum_{\tau \in \mathcal{T}} \mathbf{P}_\tau \mathbf{K}_\tau \mathbf{P}_\tau^T \mathbf{x},
}
\]
{with the local-to-global mapping $\mathbf{P}_\tau$, and the element matrices $\mathbf{K}_\tau$.}
While this already captures the main idea of matrix-free methods, further optimizations are required to reach a competitive performance.
These details are explained in great detail in the work 
\cite{KronbichlerKormann:2012a,KrKo19,MuKoKr20,CleHeiKaKr21}, and are readily included in the deal.II library \cite{dealii}.
The further optimizations include the mapping of the element-wise contributions 
{$\mathbf{K}_\tau$} %
to the reference domain.
Therein, the tensor product structure of the shape functions is exploited via sum-factorization.
Thus, even the assembly of the dense local element matrices 
{$\mathbf{K}_\tau$} %
can be replaced by efficient 
sum-factorization techniques; %
{see, e.g.,  \cite{MelenkGerdesSchwab:2001a, KronbichlerKormann:2012a}.}
This leads to a very promising computational complexity for matrix-vector multiplication, in particular for high polynomial degree shape functions.
There, we only need 
$\mathcal{O}(d^2 (p+1)(1-\frac{1}{p+1})^{-d})$ 
operations per dof compared to $\mathcal{O}((p+1)^d)$ for a sparse matrix based approach; see \cite{KronbichlerKormann:2012a}.
Furthermore, the workload is split across multiple compute nodes using MPI, and  accelerated by explicit vectorization, i.e., multiple elements are handled simultaneously with a single CPU instruction.
The operators on each level of the GMG method can be represented by such a matrix-free procedure.

The convergence rate of the multigrid method given by Algorithm~\ref{alg:mg} is 
determined
by the {\it smoothing property} and the {\it approximation property};
see \cite{Hackbusch:1985a}. In order to prove the convergence of the $W$-cycle,
it is enough to investigate the convergence of the $(\ell-(\ell-1))$ two-grid  method 
which is defined by Algorithm~\ref{alg:mg} 
{ 
when
$\mathbf{w}_{\ell - 1} \gets \MG(\mathbf{w}_{\ell - 1}, \mathbf{r}_{\ell - 1})$ 
is replaced by the exact solution
$\mathbf{w}_{\ell-1} \gets \mathbf{K}^{-1}_{\ell-1} \mathbf{r}_{\ell-1}$
of the coarse grid system $\mathbf{K}_{\ell-1} \mathbf{x}_{\ell-1} = \mathbf{r}_{\ell-1}$
in step 10.
}
In the algebraic setting, the two-grid iteration operator 
has the form
\begin{equation}
 \label{eqn:TG-IiterationMatrix}
  \mathbf{M}_\ell^{\ell-1} = \mathbf{S}_\ell^m (\mathbf{I}_\ell - \mathbf{I}_{\ell-1}^\ell \mathbf{K}_{\ell-1}^{-1} \mathbf{I}_\ell^{\ell-1}\mathbf{K}_{\ell}) \mathbf{S}_\ell^m 
  = \mathbf{S}_\ell^m \mathbf{C}_\ell^{\ell-1} \mathbf{S}_\ell^m 
 \end{equation}
where $\mathbf{I}_\ell$ is the identity matrix, 
$\mathbf{C}_\ell^{\ell-1} = \mathbf{I}_\ell - 
\mathbf{I}_{\ell-1}^\ell \mathbf{K}_{\ell-1}^{-1} \mathbf{I}_\ell^{\ell-1}\mathbf{K}_{\ell}$ 
is the coarse-grid-correction operator,
$\mathbf{S}_\ell$ denotes the iteration 
operator of the smoother, $m$ is the number of pre- and post-smoothing steps
of some symmetric smoother $\mathcal{S}_\ell$ as defined in \eqref{eqn:isus1} - \eqref{eqn:isus3}.
{Here, $\mathbf{M}_\ell^{\ell-1}$, $\mathbf{C}_\ell^{\ell-1}$, and $\mathbf{S}_\ell$ are linear operators on $X_\ell = \mathbb{R}^{N_\ell}$ with $N_\ell = \dim V_\ell + \dim Q_\ell$.}
If there is no post-smoothing, 
i.e. $\mathbf{M}_\ell^{\ell-1} =\mathbf{C}_\ell^{\ell-1} \mathbf{S}_\ell^m $,  
then two-grid convergence rate estimate
{\begin{equation}
 \label{eqn:TG-ConvergenceRate}
 \|\mathbf{M}_\ell^{\ell-1}\|_{L(X_\ell,X_\ell)} 
 \le \|\mathbf{C}_\ell^{\ell-1} \mathbf{K}_\ell^{-1}\|_{L(X_\ell^*,X_\ell)}
  \, \|\mathbf{K}_\ell \mathbf{S}_\ell^m\|_{L(X_\ell,X_\ell^*)}
 \le c \, \eta(m) =: \sigma(m) 
\end{equation}}%
follows from the {\it approximation property}
{\begin{equation}
 \label{eqn:ApproximationProperty}
 \|\mathbf{C}_\ell^{\ell-1} \mathbf{K}_\ell^{-1}\|_{L(X_\ell^*,X_\ell)} \le c,
\end{equation}}%
and the {\it smoothing property}
{\begin{equation}
 \label{eqn:SmoothingProperty}
 \|\mathbf{K}_\ell \mathbf{S}_\ell^m\|_{L(X_\ell,X_\ell^*)} \le \eta(m),
\end{equation}}%
{where 
$\|\cdot\|_{L(Y,Z)}$ denotes the standard operator norm, $\mathbf{K}_\ell$ is interpreted as a linear operator from $X_\ell$ to its dual $X_\ell^*$ and  $X_\ell = \mathbb{R}^{N_\ell}$ is equipped with a suitable mesh-dependent vector norm $\| \cdot \|_{X_\ell}$} in which we want to study the multigrid convergence,  
{$c$} is some  positive constant, and the non-negative function $\eta(m)$ converges to $0$ 
for $m$ tending to infinity. Both $c$ and $\eta$ are independent of level index $\ell$. 
Due to these properties, $\sigma(m)$ becomes smaller than $1$ for sufficiently 
many smoothing steps, and the two grid-method converges at least with the rate $\sigma(m)$.
The convergence of the $W$- and the generalized $V$-cycles then easily follows 
from a perturbation argument; see \cite{Hackbusch:1985a}.
We note that the convergence of the pure $V$-cycle is not covered by this theory.
{In the next subsection, we prove the 
approximation property \eqref{eqn:ApproximationProperty}, whereas, in Subsection~\ref{subsec:SmoothingProperty}, we propose several matrix-free smoothers,
and prove the smoothing property \eqref{eqn:SmoothingProperty}.}

\subsection{Approximation Property}
\label{subsec:ApproximationProperty}

{For the prolongation operator $\mathbf{I}_{\ell-1}^\ell$ we use the matrix representation of the canocical injection from $V_{\ell-1} \times Q_{\ell-1}$ to $V_\ell \times Q_\ell$, and set $\mathbf{I}_\ell^{\ell-1} = \left(\mathbf{I}_{\ell-1}^\ell\right)^T$ for the restriction operator.
As already pointed out, we need to specify a mesh-dependent norm $\|\cdot\|_{X_\ell}$ on $X_\ell = \mathbb{R}^{N_\ell}$. 
Let $h_\ell$ be the mesh size of $\mathcal{T}_\ell$. Then we set
\begin{equation} \label{discretenorm}
  \|\mathbf{y}_\ell\|_{X_\ell} 
    = \left( h_\ell^{d-2} \|\mathbf{v}_\ell\|^2 + h_\ell^d \|\mathbf{q}_\ell\|^2\right)^{1/2}
    \quad \text{for} \ \mathbf{y}_\ell 
    = \begin{pmatrix} \mathbf{v}_\ell \\ \mathbf{q}_\ell \end{pmatrix}
    \in X_\ell = \mathbb{R}^{N_\ell},
\end{equation}
where $\|\cdot\|$ denotes the Euclidean norm of vectors.
If $\mathbf{v}_\ell$ and $\mathbf{q}_\ell$ are the vector representations of some elements $v_\ell \in V_\ell$ and $q_\ell \in Q_\ell$, respectively, then it is easy to see that $\|\mathbf{y}_\ell\|_{X_\ell}$ and $\|(v_\ell,q_\ell)\|_{0,\ell}$, given by
\[
  \|(v_\ell,q_\ell)\|_{0,\ell} 
    = \left( h_\ell^{-2} \|v_\ell\|_{L^2(\Omega)^d}^2 + \|q_\ell\|_{L^2(\Omega)}^2 \|\right)^{1/2}, \quad \text{for} \, (v_\ell,q_\ell) \in V_\ell \times Q_\ell,
\]
are equivalent norms, in short
\[
  \|\mathbf{y}_\ell\|_{X_\ell} \sim \|(v_\ell,q_\ell)\|_{0,\ell}.
\]
Associated with this norm on $V_\ell\times Q_\ell$ and the bilinear form $\mathcal{B}$,
a second mesh-dependent norm  on $V_\ell\times Q_\ell$ is introduced by
\begin{equation*} %
  \|(v_\ell,q_\ell)\|_{2,\ell} = \sup_{0 \neq (w_\ell,r_\ell) \in V_\ell \times Q_\ell} \frac{\mathcal{B}((v_\ell,q_\ell),(w_\ell,r_\ell))}{\|(w_\ell,r_\ell)\|_{0,\ell}}
  \quad \text{for} \, (v_\ell,q_\ell) \in V_\ell \times Q_\ell.
\end{equation*}
Obviously, $\|(v_\ell,q_\ell)\|_{2,\ell}$ and $\|\mathbf{K}_\ell \mathbf{y}_\ell\|_{X_\ell^*}$ are equivalent norms, too, i.e.,
\[
  \|(v_\ell,q_\ell)\|_{2,\ell} \sim \|\mathbf{K}_\ell \mathbf{y}_\ell\|_{X_\ell^*}, 
\]
where $\|\cdot\|_{X_\ell^*}$ denotes the norm in $X_\ell^*$ dual to $\|\cdot\|_{X_\ell}$.

Observe that $\mathbf{K}_{\ell-1}^{-1} \mathbf{I}_\ell^{\ell-1}\mathbf{K}_{\ell}$ is the matrix representation of the Ritz projection $R_\ell^{\ell-1}$ from $V_\ell \times Q_\ell$ onto $V_{\ell-1}  \times Q_{\ell-1}$, given by
\[
  \mathcal{B}(R_{\ell}^{\ell-1}(v_\ell,q_\ell),(w_{\ell-1},r_{\ell-1})) = \mathcal{B}((v_\ell,q_\ell),(w_{\ell-1},r_{\ell-1})),
\]
for all $(v_\ell,q_\ell) \in V_\ell \times Q_\ell$ and $(w_{\ell-1},r_{\ell-1}) \in V_{\ell-1} \times Q_{\ell-1}$.
Under the assumption of elliptic regularity the  approximation property was shown in \cite{Verfuerth1984} in the following form.
\[
  \|[I - R_{\ell}^{\ell-1}] (v_\ell,q_\ell)\|_{0,\ell} \le c \, \|(v_\ell,q_\ell)\|_{2,\ell} \quad \text{for all} \  (v_\ell,q_\ell) \in V_\ell \times Q_\ell,
\]
with a constant $c$ which is independent of the level $\ell$. Therefore,
\[
  \|\mathbf{C}_\ell^{\ell-1} \mathbf{y}_\ell \|_{X_\ell}
   \sim \|[I - R_{\ell}^{\ell-1}] (v_\ell,q_\ell)\|_{0,\ell}
   \le c \, \|(v_\ell,q_\ell)\|_{2,\ell} 
   \sim c \, \|\mathbf{K}_\ell \mathbf{y}_\ell\|_{X_\ell^*}
   \quad \text{for all} \ \mathbf{y}_\ell \in X_\ell,
\]
which immediately implies the approximation property of the form \eqref{eqn:ApproximationProperty}.
}

\subsection{Matrix-free Smoothers and  Smoothing Property}
\label{subsec:SmoothingProperty}

We now introduce different classes of properly scaled Chebyshev-Jacobi smoothers,
which allow an efficient matrix-free implementation, 
and investigate their smoothing property. 
Since the smoothers are living only on one level, we omit the level subindex  $\ell$ 
throughout this subsection.

\subsubsection{Scaled Chebyshev-Jacobi Smoothers}
\label{subsubsec:ScaledChebyshevJacobiSmoothers}

The main challenge of a matrix-free implementation  of the GMG algorithm is connected 
with an appropriate smoother that must permit an efficient matrix-free implementation as well.
Since the matrix entries are not at our disposal, many classical methods like 
{(block)-Gauss-Seidel, (block)-Jacobi,} or ILU, are not applicable.
Rather, we 
should look for smoothing 
methods that only require matrix-vector operations.
One class of methods that seems to be predestined for this task are polynomial methods.
{This approach yields}
smoothers that are based on a matrix polynomial 
$s_{k_C}(C)$ of the degree $k_C$ .
Such a polynomial can be efficiently applied by only performing matrix-vector multiplications with the matrix-free operator  
$C$.

We can further improve the situation by precomputing the diagonal of a matrix-free operator, which enables us to use Jacobi-based methods.
The diagonal can be computed efficiently using the deal.II matrix-free framework, and storing it comes at the cost of only one additional vector.
{ 
In fact, we will use the Jacobi-method, and improve its smoothing property by means of 
polynomial acceleration.
}

For a matrix-free implementation, our choices for $\hat{A}$ and $\hat{S}$ 
in \eqref{eqn:isus1}--\eqref{eqn:isus3}
are focused on Chebyshev-Jacobi smoothers for the realization of 
{matrix-free application of
$\hat{A}^{-1}$ and $\hat{S}^{-1}$
to vectors.}

For a linear system of the form $M x = b$ with some (simple) preconditioner $D$ these smoothers are semi-iterative algorithms of the form 
\[
	\tilde{x} = x + s_{k}(D^{-1} M) D^{-1} (b - Mx)
\]
with 
$s_{k}(t) = (1 - T_{k+1}(t)) t^{-1}$ and 
$T_{k+1}(t) := {C_{k+1}(\frac{\beta + \alpha - 2t}{\beta - \alpha})}/{C_{k+1}(\frac{\beta + \alpha}{\beta - \alpha})},$
where $C_{k+1}$ denotes the Chebyshev polynomial of the first kind of degree $k+1$. %
We note that, in this subsection, the subindex $k$ is related to the degree 
of the Chebyshev polynomial and not to the polynomial degree of the shape functions
used for constructing the finite element spaces $\mathcal{Q}_{k}$-$\mathcal{Q}_{k-1}$ introduced in 
Section\ref{sec:mVF+mFEM}. Later we will again add the subindex that indicates 
the matrix to which the polynom belongs as we did above with $k_C$.
The interval $[\alpha, \beta]$ contains the high energy eigenvalues of $D^{-1} M$. 
In particular, we will use $[\alpha,\beta] = [\beta/2,\beta]$ with $\beta = \lambda_\text{max}(D^{-1} M)$. 
Actually, in the numerical tests, a close approximation of $\lambda_\text{max}(D^{-1} M)$ will be used for $\beta$.
The iteration matrix for the Chebyshev-Jacobi method is then given by $T_{k+1}(D^{-1} M)$ with associated preconditioner $C_M^{-1} = \Cheb(M,D,k) = s_{k}(D^{-1} M) D^{-1}$.
The Chebyshev-Jacobi method is easy to apply by exploiting a three-term recurrence relation,
see, e.g., \cite{Saad:2003a}.
For $k=0$, it corresponds to the damped Jacobi method with damping parameter 
$\omega = 2/(\alpha + \beta)$
{ 
and the iteration matrix $T_{1}(D^{-1} M) = I - \omega D^{-1}M$.
}

The next lemma contains two important estimates for the Chebyshev-Jacobi smoother, which are required for the scaling of the smoother and the proof of the smoothing property.

\begin{Lemma} \label{lemmaone}
Let $M$ and $D$ be symmetric and positive definite matrices, $k \in \mathbb{N}_0$, and $\alpha, \beta \in \mathbb{R}$ with $0 < \alpha < \beta = \lambda_\text{max}(D^{-1} M)$. Then the following estimates hold.
\[
  c_1 \, M \le C_M 
  \quad \text{and} \quad
  C_M \le c_2 \, \beta \, D,
\]
where
\[
  c_1 = \left(1 + \frac{1}{C_{k+1}(\bar s)}\right)^{-1}, \quad
  c_2 = \left( 1 -  \frac{1}{C_{k+1}(\bar s)}\right)^{-1}
\]
with $\bar s = (\beta+\alpha)/(\beta-\alpha)$.
\end{Lemma}
\begin{proof}
We have
\begin{align*}
  M^{1/2} C_M^{-1} M^{1/2}
   & = M^{1/2} s_{k}(D^{-1} M) D^{-1} M^{1/2}
     =  s_{k}(M^{1/2} D^{-1} M^{1/2}) M^{1/2} D^{-1} M^{1/2} \\
   & \le \lambda_\text{max} \left( s_{k}(M^{1/2} D^{-1} M^{1/2}) M^{1/2} D^{-1} M^{1/2} \right) \, I
\end{align*}
Since $\sigma(M^{1/2} D^{-1} M^{1/2}) = \sigma(D^{-1} M) \subset [0,\beta]$, 
we immediately obtain
\begin{equation*}
	\lambda_\text{max} \left( s_{k}(M^{1/2} D^{-1} M^{1/2}) M^{1/2} D^{-1} M^{1/2} \right)
	\leq \max_{\lambda \in [0,\beta]} \left( s_k(\lambda) \, \lambda \right) =
	  1 + \frac{1}{C_{k+1}(\bar s)}  = \frac{1}{c_1},
\end{equation*}
which proves the first inequality.

For the second inequality,  we proceed in a similar way. We obtain
\begin{align*}
  D^{1/2} C_M^{-1} D^{1/2}
    & = D^{1/2} s_{k}(D^{-1} M) D^{-1} D^{1/2}
      =  s_{k}(D^{-1/2} M D^{-1/2}) \\
    & \ge \lambda_\text{min}\left(s_{k}(D^{-1/2} M D^{-1/2})\right) \, I
\end{align*}
and
\[
  \lambda_\text{min}\left(s_{k}(D^{-1/2} M D^{-1/2})\right) 
    \ge \min_{\lambda \in [0,\beta]} s_{k}(\lambda) 
    \ge \frac{1}{\beta} \left( 1 -  \frac{1}{C_{k+1}(\bar s)}\right) = \frac{1}{c_2 \beta}, 
\]
which completes the proof of the second inequality. 
The used bounds involving the polynomial $s_k$ follow easily from well-known properties of Chebyshev polynomials.
\end{proof}

Note that, for a fixed ratio $\alpha/\beta$, the quantities $c_1$ %
and $c_2$ in Lemma \ref{lemmaone} depend only on $k$ %
and %
approach 1 for $k$ approaching $\infty$. %

With these notations we set
\[
  \hat{A}^{-1} = \sigma \, \Cheb(A,D_A,k_A)
\]
with $D_A := \diag(A)$, $\beta =\lambda_\text{max}(D_A^{-1}A)$, $\alpha = \beta/2$, and $\sigma = c_1 = c_1(k_A,\alpha,\beta)$ according to Lemma \ref{lemmaone}, which guarantees $\hat{A} \geq A$.

In principle the same construction could be carried out for $\hat{S}$ leading to
\[
  \hat{S}^{-1} = \tau \, \Cheb(\tilde{S}, D_{\tilde S}, k_S)
\]
with $\tilde{S}= B\hat{A}^{-1}B^T$, $D_{\tilde S} := \diag(\tilde{S})$, $\beta = \lambda_\text{max}(D_{\tilde S}^{-1} \tilde S)$, $\alpha = \beta/2$,  and $\tau = c_1 = c_1(k_S,\alpha,\beta)$ according to Lemma \ref{lemmaone}, which guarantees $\hat{S} \ge \tilde{S}$.
However, obtaining $\diag(\tilde{S})$ as required for the Chebyshev-Jacobi smoother is more expensive, 
since it involves the inverse $\hat{A}^{-1}$.
Hence, we aim to replace $D_{\tilde S}$ %
by cheaper approximations, namely
\begin{align}
	\tilde D_{\tilde S_d} &:= \diag(\tilde{S}_d) 
	\quad \text{with} \quad \tilde{S}_d := B (\diag(A))^{-1} B^T, \label{eq:diag d} \\
	\tilde D_{\tilde S_p} &:= \diag(\tilde{S}_p) %
	\quad \text{with} \quad \tilde{S}_p := M_p,
	\label{eq:diag p}\\
	\tilde D_{\tilde S_{loc}} &:= \sum_{\tau \in \mathcal{T}_h} P_\tau^T 
	{\diag( B_\tau (\diag(A_\tau))^{-1} B_\tau^T ) P_\tau,\label{eq:diag loc}}
\end{align}
where $M_p$ denotes the pressure mass matrix, 
and where the local-to-global matrices $P_\tau$ were defined before.
The approximation $\tilde D_{\tilde S_{loc}}$ computes 
the diagonal
$\diag( B_\tau (\diag(A_\tau))^{-1} B_\tau^T )$ locally on each element $\tau$ and assembles it together.
The computation of $\tilde D_{\tilde S_d}$ still requires to assemble the sparse matrix $B$, in order to compute $(\tilde{D}_{\tilde S_d})_{ii} = \sum_{j} B_{ij}^2 (\diag(A)_{jj})^{-1}$.
Thus, in order to keep the smoother fully matrix-free, the choices $\tilde D_{\tilde S_p}$ and $\tilde D_{\tilde S_{loc}}$ are preferred.

{
The setup of the ingredients of the smoother is summarized in Algorithm~\ref{alg:initialization}.	
The smoother \eqref{eqn:isus1} -  \eqref{eqn:isus3} with $\hat{A}$ and $\hat{S}$ 
defined by Algorithm~\ref{alg:initialization} can be written in the compact form
\begin{equation}
\label{eqn:SmootherCompactForm}
 \mathbf{x}^{j+1} = \mathbf{x}^{j} + \mathbf{\hat{K}}^{-1} [\mathbf{b} - \mathbf{K}  \mathbf{x}^{j} ]
\end{equation}
with 
\begin{equation*}
\label{eqn:hatK}
\mathbf{\hat{K}} = \begin{pmatrix} \hat{A} & B^T\\ B & B\hat{A}^{-1}B^T - \hat{S}  \end{pmatrix}
\end{equation*}
The iteration matrix $\mathbf{S}$ of the smoother \eqref{eqn:SmootherCompactForm} 
is obviously given by $ \mathbf{S} = \mathbf{I} - \mathbf{\hat{K}}^{-1}\mathbf{K}$.
We note that, in Algorithm~\ref{alg:initialization}, 
$\bar s = 3$ for the standard  choice $\alpha= \beta/2$.
}	
\begin{algorithm}
		\SetKwFunction{MG}{MG}
		\SetKwInOut{Input}{Input}	
		
		\Input{$k_A, k_S$}
		
		\Indp
		
		\vspace{2mm}
		\tcc{Velocity part}
		
		$D_A \gets \diag(A)$

		$\beta \gets \lambda_{max}(D_A^{-1} A)$
		
		$\alpha \gets \beta/2$
		
		$C_A^{-1} \gets \Cheb(A, D_A, k_A)$
		
		$\sigma = \frac{C_{k_A+1}(\bar s)}{1 + C_{k_A+1}(\bar s)}$
		
		$\hat{A}^{-1} \gets \sigma \, C_A^{-1}$\tcc*{scale Chebyshev to guarantee $\hat{A} \geq A$}
		
		\vspace{2mm}
		\tcc{Pressure part}
		
		$\tilde D_{\tilde{S}} \gets \tilde D_{\tilde S_d}, \tilde D_{\tilde S_p},$ or $\tilde D_{\tilde S_{loc}}$\tcc*{see \eqref{eq:diag d},\eqref{eq:diag p}, or \eqref{eq:diag loc}}

		$\beta \gets \rho(\tilde D_{\tilde{S}}^{-1} \tilde{S})$
		
		$\alpha = \beta/2$
		
		$C_S^{-1} \gets Cheb(\tilde{S},\tilde D_{\tilde{S}}, k_S)$
		
		$\tau = \frac{C_{k_S+1}(\bar s)}{1 + C_{k_S+1}(\bar s)}$
		
		$\hat{S}^{-1} \gets \tau \,  C_S^{-1}$\tcc*{scale
		Chebyshev to guarantee $\hat{S} \geq \tilde{S}$}
		
		\vspace{2mm}
		
		\Return{$\hat{A}^{-1}, \hat{S}^{-1}$}
		
		\caption{Setup of the smoother for the Stokes problem.}
		\label{alg:initialization}
	\end{algorithm}

\subsubsection{Smoothing Property}
\label{subsubsec:SmoothingProperty}

{
Since we can always ensure by scaling that the ingredients $\hat{A}$ and $\hat{S}$ of 
$\mathbf{\hat{K}}$ satisfy the spectral inequalities \eqref{eqn:Smoother:SpectralInequalities}
for all choices proposed in Subsection~\ref{subsubsec:ScaledChebyshevJacobiSmoothers},
Theorem~7 from \cite{DrzisgaJohnRuedeWohlmuthZulehner:2018a} 
(see also Theorem~3 in \cite{SchoeberlZulehner:2003a})
delivers the estimate
\begin{equation}
 \label{eqn:SmoothingPreProperty}
 \|\mathbf{K} \mathbf{S}^m\|_{L(X_\ell,X_\ell^*)} \le \eta_0(m-1)\, \| \mathbf{\hat{K}} - \mathbf{K} \|_{L(X_\ell,X_\ell^*)}
\end{equation}
with
\[
  \eta_0(m) = \frac{1}{2^m} \binom{m}{ \lfloor (m+1)/2 \rfloor}
    \le \left\{ \begin{array}{ll}
                  \sqrt{\displaystyle
		        \frac{2}{\pi m}}  & \text{for even }  m ,\\   %
                  \sqrt{\displaystyle
		        \frac{2}{\pi (m+1)}}  & \text{for odd } m ,   %
	       \end{array} \right. 
\]
where $\binom{m}{\ell}$ denotes the binomial coefficient, 
and $\lfloor x \rfloor$ 
is nothing but
the largest integer less or equal to $x \in \mathbb{R}$.

It remains to show that
\begin{equation}
\label{eqn:K-hatK}
\| \mathbf{\hat{K}} - \mathbf{K} \|_{L(X_\ell,X_\ell^*)} \le c 
\end{equation}
to get the smoothing property \eqref{eqn:SmoothingProperty} 
with $\eta(m) = c \eta_0(m) = O(1/\sqrt{m})$,
where $c$ is some non-negative constant 
that should be independent of the level index $\ell$. 
Indeed, following \cite{SchoeberlZulehner:2003a}, we can proceed 
as follows:} 
First we note that
\[
  \mathbf{\hat{K}} - \mathbf{K} 
    = \begin{pmatrix}
         \hat{A} - A & 0 \\ 0 & \hat{S} - \tilde{S}
      \end{pmatrix}.
\]
Therefore, it follows from the definition \eqref{discretenorm} of the norm in $X_\ell$ that
\begin{align*}
  \| \mathbf{\hat{K}} - \mathbf{K} \|_{L(X_\ell,X_\ell^*)}
    & = \max \left( h^{-(d-2)} \|\hat{A}-A\|, h^{-d} \|\hat{S} - \tilde{S}\| \right) \\
    & \le \max \left( h^{-(d-2)} \|\hat{A}\|, h^{-d} \|\hat{S}\| \right),
\end{align*}
where here the matrix norm $\| \cdot \|$ is the spectral norm that is 
associated to the Euclidean vector norm. So, \eqref{eqn:K-hatK} holds, if
\[
  h^{-(d-2)} \|\hat{A}\| \le c
  \quad \text{and} \quad
  h^{-d} \|\hat{S}\|.
\]
We start with the verification of the first estimate
\begin{equation} \label{firstest}
  h^{-(d-2)} \|\hat{A}\| \le c.
\end{equation}
From Lemma \ref{lemmaone} 
with $M = \hat{A}$ and $D = D_A$,
it follows that
\[
  \hat{A} \le \frac{c_2}{c_1} \, \lambda_\text{max}(D_A^{-1} A) \, D_A,
\]
which leads to $\|\hat{A}\| \le (c_2/c_1) \, \lambda_\text{max}(D_A^{-1} A) \|D_A\|$. In order to conclude from this estimate that \eqref{firstest} holds
for some constant $c$ independent of the mesh level, it suffices to show that
\[
  A \le c \, D_A \quad \text{and} \quad \|D_A\| \le c \, h^{d-2}.
\]
The first part, which is equivalent to $\lambda_\text{max}(D_A^{-1} A) \le c$, follows, e.g., from \cite[Theorem 12.20] {Hackbusch2016} by exploiting the sparsity pattern on $A$,
whereas the second part can easily be  obtained by a standard scaling argument.

Next we discuss the second estimate
\[
  h^{-d} \|\hat{S}\| \le c
\]
with $\hat{S}^{-1} = \tau \, \Cheb(\tilde S, \tilde D_S,k_S)$ and the three different choices for $\tilde D_S$ 
given by \eqref{eq:diag d}, \eqref{eq:diag p}, and \eqref{eq:diag loc}.
For each case it follows from Lemma \ref{lemmaone} that $\|\hat{S}\| \le (c_2/c_1) \, \lambda_\text{max}(\tilde D_S^{-1} \tilde{S}) \|\tilde D_S\|$. As before it suffices to show that
\begin{equation} \label{Ssufficient}
  \tilde S  \le c \, \tilde D_S \quad \text{and} \quad
  \|\tilde D_S\| \le c \, h^d.
\end{equation}
We will discuss these conditions for each of the three choices for $\tilde D_S$ separately.
\begin{enumerate}
\item
We start with the choice $\tilde D_S = \tilde D_{\tilde S_p} = \diag(M_p)$. 
The coercivity of the bilinear form $a(\cdot,\cdot)$ and the boundedness of the bilinear form $b(\cdot,\cdot)$ 
yield the estimate
\[
  \sup_{v_h \in V_h} \frac{b(v_h,q_h)^2}{a(v_h,v_h)} \le c \, \|q_h\|_{L^2(\Omega)}^2,
\]
or, equivalently, in matrix notation,
\[
  B A^{-1} B^T \le c \, M_p. 
\]
From the scaling condition $\hat{A} \ge A$, we obtain the following upper bound of $\tilde S$.
\[
  \tilde S = B \hat{A}^{-1} B^T \le  B A^{-1} B^T \le c \, M_p. 
\]
As before by exploting the sparsity pattern of $M_p$ and by using a standard scaling argument it follows that
\[
  M_p \le c \, \diag(M_p) \quad \text{and} \quad \|\diag(M_p)\| \le c h^d,
\]
which, together with $\tilde S \le c \, M_p$, directly yield the required estimates \eqref{Ssufficient}.
\item
Next we discuss the choice $\tilde D_S = \tilde D_{\tilde S_d} = \diag(B (\diag(A))^{-1} B^T)$. 
It suffices to show that 
\[
  \tilde D_{\tilde S_d} \sim \tilde D_{\tilde S_p}.
\]
Then the estimates \eqref{Ssufficient} follow directly from the previous case $\tilde D_S = \tilde D_{\tilde S_p}$.

For $v_h \in V_h$ and $q_h \in Q_h$, we have 
\[
  b(v_h,q_h) = - \int_\Omega q_h \div v_h \ dx = \int_\Omega v_h \cdot \nabla q_h \ d x 
\]
via
integration by parts, which in turn implies
\[
  \sup_{v_h \in V_h} \frac{b(v_h,q_h)}{\|v_h\|_{L^2(\Omega)^d}} \le \|\nabla q_h\|_{L^2(\Omega)},
\]
or, equivalently, in matrix notation
\[
  B M_v^{-1} B^T \le K_p,
\]
where $M_v$ is the mass matrix in $V_h$ and $K_p$ is the matrix representing the $H^1$-seminorm in $Q_h$. A reverse inequality of the form
\[
  B M_v^{-1} B^T \gtrsim \, K_p
\]
also holds, see \eqref{BerPirlocalMatrixVersion}. Therefore
we have
\[
  B M_v^{-1} B^T \sim K_p.
\]
Here we used the following notation for symmetric matrices $M$, $N$: $M \lesssim N$ if there is a constant $c > 0$ such that $M \le c \, N$, $M \gtrsim N$ if $N \lesssim M$, and $M \sim N$ if $M \lesssim N$ and $M \gtrsim N$. The involved constants are meant to be independent of the mesh level.

A standard scaling argument gives
\[
  \diag(A) \sim h^{-2} \, M_v.
\]
Hence
\[
  B (\diag(A))^{-1} B^T \sim h^{2} B M_v^{-1} B^T \sim h^{2}\, K_p
\]
and
\[
  \tilde D_S = \diag(B (\diag(A))^{-1} B^T) \sim h^{2} \, \diag(K_p) \sim \diag(M_p) = \tilde D_{\tilde S_p},
\]
which concludes the discussion of this case.
\item
For the choice $\tilde D_S = \tilde D_{\tilde S_{loc}} = \sum_{\tau \in \mathcal{T}} P_\tau^T \diag( B_\tau (\diag(A_\tau))^{-1} B_\tau^T ) P_\tau$ we procede similarly to the previous case:
We have, see \eqref{BerPirlocalMatrixVersion},
\[
  B_\tau M_{v,\tau}^{-1} B_\tau^\top
   \gtrsim K_{p,\tau}.
\]
An inequality in reverse direction
\[
   B_\tau M_{v,\tau}^{-1} B_\tau^\top 
   \lesssim  h_\tau^{-2}  M_{p,\tau}
\]
follows from a standard scaling argument. These two estimates imply 
\[
   h_\tau^{-2}  \diag(M_{p,\tau})
   \sim \diag(K_{p,\tau})
   \lesssim \diag \left( B_\tau M_{v,\tau}^{-1} B_\tau^\top \right)
   \lesssim h_\tau^{-2}  \diag(M_{p,\tau}).
\]
Therefore,
\[
   \diag \left( B_\tau M_{v,\tau}^{-1} B_\tau^\top \right)
   \sim h_\tau^{-2}  \diag(M_{p,\tau})
\]
and
\begin{align*}
  & \sum_{\tau \in \mathcal{T}_h} 
    P_\tau^T  \diag( B_\tau (\diag(A_\tau))^{-1} B_\tau^T ) P_\tau \\
  & \quad \sim \sum_{\tau \in \mathcal{T}_h} 
    h^2 \, P_\tau^T  \diag( B_\tau M_{v,\tau}^{-1} B_\tau^T ) P_\tau 
    \sim \sum_{\tau \in \mathcal{T}_h} 
    P_\tau^T  \diag(M_{p,\tau}) P_\tau \sim \diag(M_p).
\end{align*}
We have shown that
\[
  \tilde D_{\tilde S_{loc}} 
  = \sum_{\tau \in \mathcal{T}_h} 
    P_\tau^T  \diag( B_\tau (\diag(A_\tau))^{-1} B_\tau^T ) P_\tau
  \sim \diag(M_p) = \tilde D_{\tilde S_p},
\]
which concludes the discussion of the last case.
\end{enumerate}

\section{Numerical Results}
\label{sec:NumericalResults}
In this section, we 
perform
some numerical experiments in two and three spatial dimensions 
in order to 
illustrate
and quantify the theoretical results.
First, the numerical setting is described. Then we present the results of
extensive numerical tests for 
the unit square with homogeneous Dirichlet boundary conditions 
and some right-hand side force,
typical benchmark problems, namely the driven cavity Stokes problem in 2d 
and 3d,  the 2d Stokes flow around a cylinder, and finally the 2d instationary driven cavity Stokes problem
with implicit time discretization of which leads to the solution of a generalized 
Stokes problem at every time step.

\subsection{Numerical Setting and Programming Code}
\label{subsec:NumericalSetting+ProgrammingCode}
The code is implemented in C++ using the open-source finite element library deal.II \cite{dealii}, and in particular, their matrix-free and geometric multigrid components; see 
\cite{CleHeiKaKr21,KronbichlerKormann:2012a}.
For the GMG method, we use a W-cycle with an equal number of pre- and post-smoothing steps, which will be varied in the numerical experiments.
On the coarsest grid, we assemble the sparse matrix and apply a direct method 
for solving the algebraic system.
This is feasible due to the small amount of dofs on $\mathcal{T}_{\ell = 0}$.
For more complex geometries 
yielding more dofs on the coarse grid $\mathcal{T}_{0}$,
iterative solvers might be a more suitable and more efficient alternative.
In all computations presented below,
the finite element discretization 
is always done by 
inf-sup stable $\mathcal{Q}_2$-$\mathcal{Q}_1$ elements
on quadrilateral and hexahedral meshes in two and three spatial dimensions, respectively.

In order to measure the multigrid convergence rate $\sigma(m)$,
we perform $i \leq 15$ iterations of the multigrid solver until a reduction of the initial residual by $10^{-10}$ is achieved.
Since the first few steps usually yield a much better reduction compared to the later iterations, we only consider the average reduction $q$ of the last $3$ steps; see, e.g., \cite{TrOoSc01}.
Hence, we compute
\[
	q = \left( \frac{r^j}{r^{j-3}} \right)^{\frac{1}{3}},
\]
where $r^j = \| \mathbf{b} - \mathbf{Kx}^j\|$ denotes the residual after $j$ GMG iterations, 
with the initial guess $\mathbf{x}^0 = \mathbf{0}$.

For the timings, we report the time $T$ needed to reduce the error by $\varepsilon = 10^{-1}$, i.e.,
\[
	T(\varepsilon) = T_{iter} \frac{\log(\varepsilon)}{\log(q)},
\]
with the time for a single iteration $T_{iter}$ and the respective reduction rate $q$ obtained as described before.
The experiments were done on one node of our local hpc cluster Radon1\footnote{https://www.oeaw.ac.at/ricam/hpc}, 
equipped with a Intel(R) Xeon(R) Silver 4214R CPU with 2.40GHz and 128 GB of RAM.

\subsection{Laser beam material processing}
\label{subsec:UnitSquarePointForce}
In order to substantiate our theoretical setting, we first consider 
an example with homogeneous Dirichlet boundary conditions. 
The background of this example is in laser material processing; see 
for instance \cite{OTTO201035}. Therein, various physics interact in which the heated 
material can be modeled with the help of viscous flow (here Stokes). 
Concentrating on the fluid component only, the laser beam is mathematically  modeled 
by a non-homogeneous right-hand side point force that is here represented by a smoothed 
Gaussian curve.

\subsubsection{Geometry, Boundary Conditions, Parameters, and Right Hand Side}
\begin{figure}[htb]
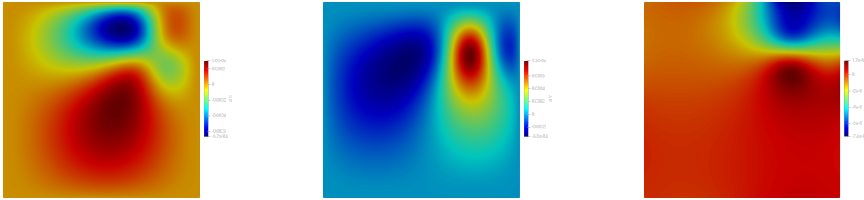

	\centering
	\includegraphics[width=0.2\textwidth]{img/laser_x}
	\hspace{1cm}
	\includegraphics[width=0.2\textwidth]{img/laser_y}
	\hspace{1cm}
	\includegraphics[width=0.2\textwidth]{img/laser_p}
	\caption{Laser beam material processing: Computed solution of the first test case. (a) x-velocity, (b) y-velocity, (c) pressure}
	\label{img:laser}
\end{figure}

The domain is the unit square $(0,1)^2$. For the boundary 
conditions we have $u_D =  (0, 0)^T$.
The (dynamic) viscosity is given by $\mu = 0.0022$.
Given the point $x_0=(0.75,0.75)\in \Omega$, the right hand side $f$
of \eqref{eqn:Stokes:cf} is given by
\begin{equation}
f = 
\begin{pmatrix}
f_1 \\
f_2
\end{pmatrix}
=
\begin{pmatrix}
0 \\
c_1 \exp(-c_2 (x-x_0)^2)
\end{pmatrix}
, \quad x\in\Omega,
\end{equation}
with $c_1 = 0.001$ and $c_2 = 100$.
Figure \ref{img:laser} illustrates the computed solutions.

\subsubsection{Results for Different Diagonal Approximations}

\begin{figure}[htb]
	\centering
	\includegraphics[width=0.45\textwidth]{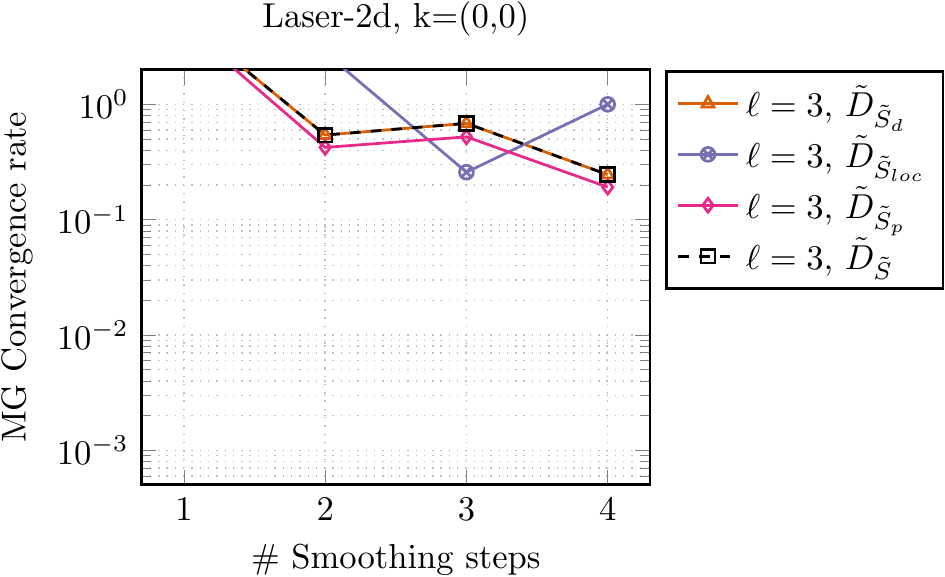}
	\includegraphics[width=0.45\textwidth]{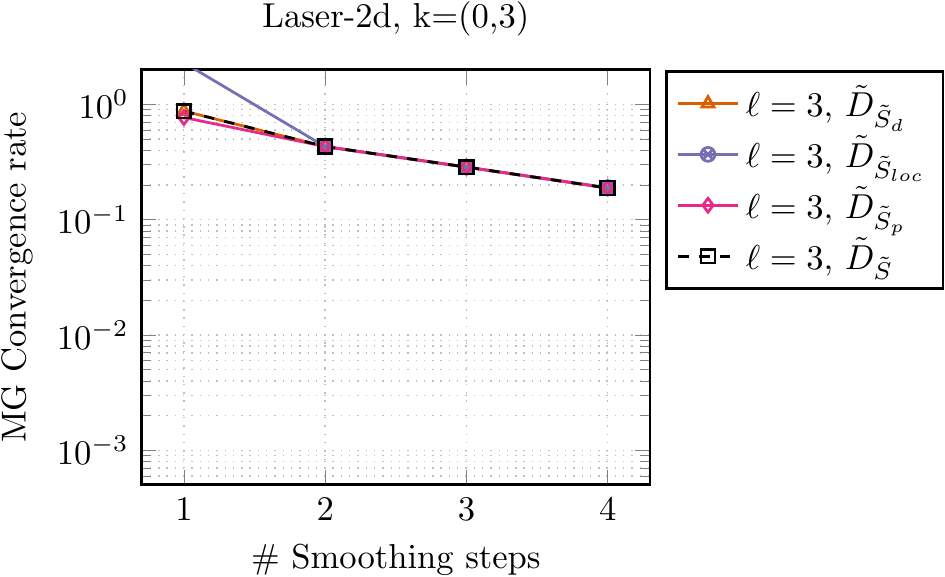}
	
	\vspace{1mm}
	
	\includegraphics[width=0.45\textwidth]{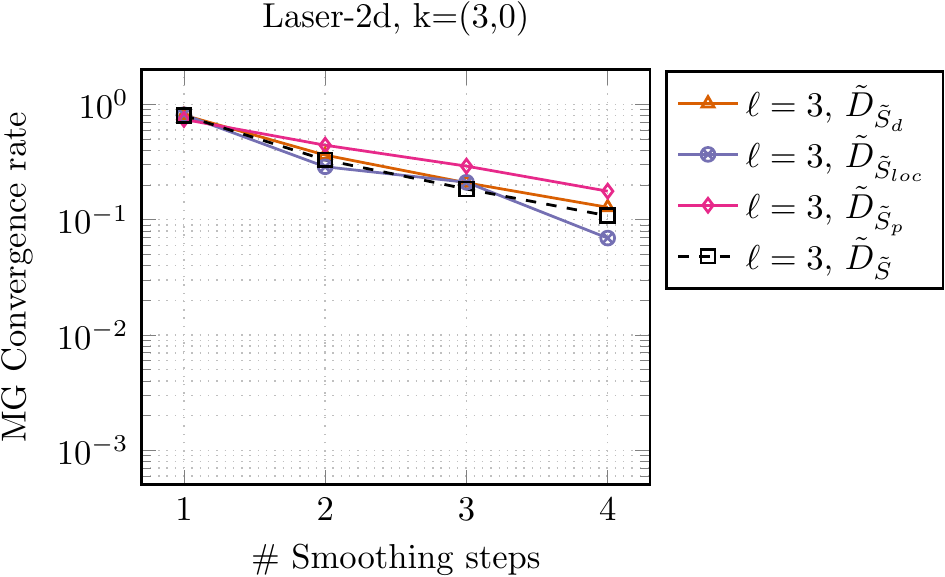}
	\includegraphics[width=0.45\textwidth]{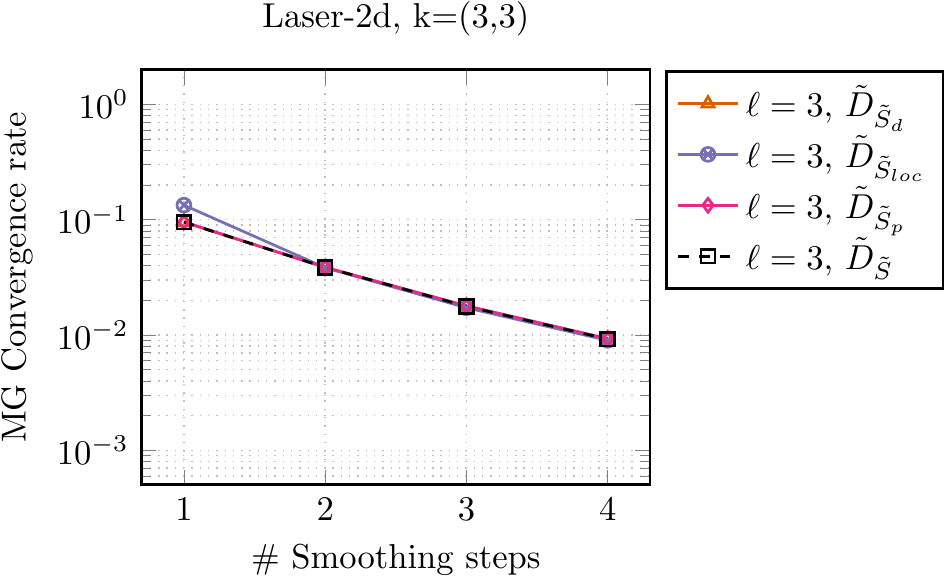}
	\caption{Laser beam material processing: Convergence rates for different diagonal approximations of the Schur complement
		and various Chebyshev polynomial degrees.}
	\label{fig:laser:diagonal}
\end{figure}

We start with comparing the convergence rates using the different approximations for $\tilde{D}_{\tilde{S}}$ 
which are visualized in Figure \ref{fig:laser:diagonal}.
For the $A$-block, we can easily compute $D_A = \diag(A)$ and, hence, apply $\hat{A}^{-1} := \hat\sigma C_A^{-1}$.
However, for the Schur-complement part, computing the exact diagonal $D_{\tilde{S}}$ is 
computationally expensive
since it involves $\hat{A}^{-1}$.
Hence, in 
\eqref{eq:diag d}, \eqref{eq:diag p}, and \eqref{eq:diag loc},
we suggested several approximations $\tilde{D}_S$, which are visualized in 
Figure \ref{fig:laser:diagonal}.
We here consider only $2467$ dofs since we want to compare different cheaper approximations 
with the expensive diagonal $D_{\tilde{S}}$.
The plot only shows the results for the Jacobi setting ($k=0$) and the Chebyshev-Jacobi degrees $k=3$, but the behavior is similar for other polynomial degrees in between,
and will be investigated later on in the more complex examples.
We can see that using any of the approximations $\tilde{D}_S$ yields mostly similar rates as the exact diagonal $D_{\tilde S}$.
In any case, increasing the number of smoothing steps improves the convergence rate of the multigrid solver.
In the Jacobi case 
presented in top-left corner of Figure~\ref{fig:laser:diagonal},
this is not obvious from the plot.
However, for sufficiently many smoothing steps, the expected behavior $q \sim 1/\sqrt{m}$ is attained.

\clearpage
\subsection{2d Driven Cavity}
\label{subsec:DrivenCavity2d}
The second example deals with the stationary Lid-Driven Cavity flow \cite{cavity:review,ZiTaNi06};
see also the benchmark collection of TU Dortmund.\footnote{\url{http://www.mathematik.tu-dortmund.de/~featflow/album/catalog/ldc_low_2d/data.html}}

\subsubsection{Geometry, Boundary conditions, and Parameters}
\label{subsubsec:DrivenCavity2d:Data}
\begin{figure}[htb]
	\centering
	\includegraphics[width=0.2\textwidth]{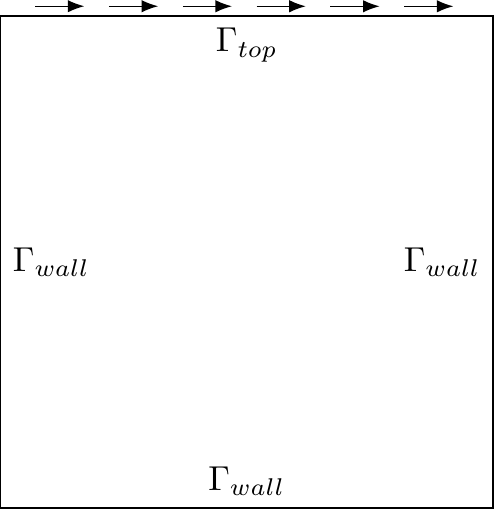}
	\hspace{1cm}
	\includegraphics[width=0.2\textwidth]{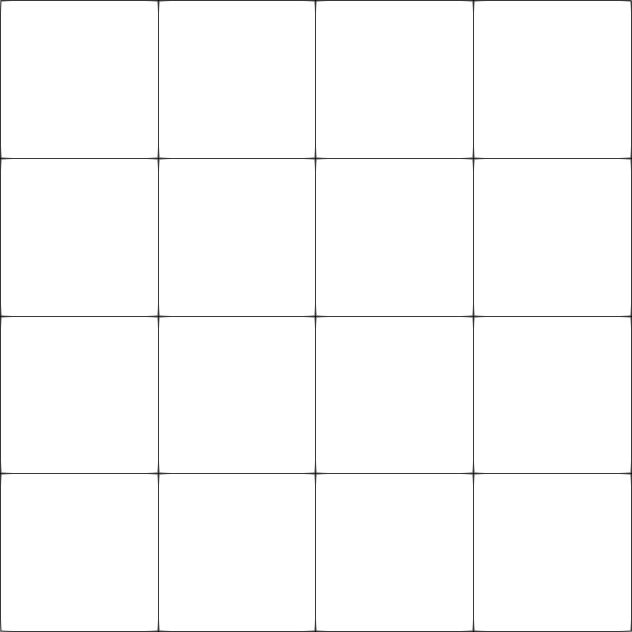}
	\hspace{1cm}
	\includegraphics[width=0.2\textwidth]{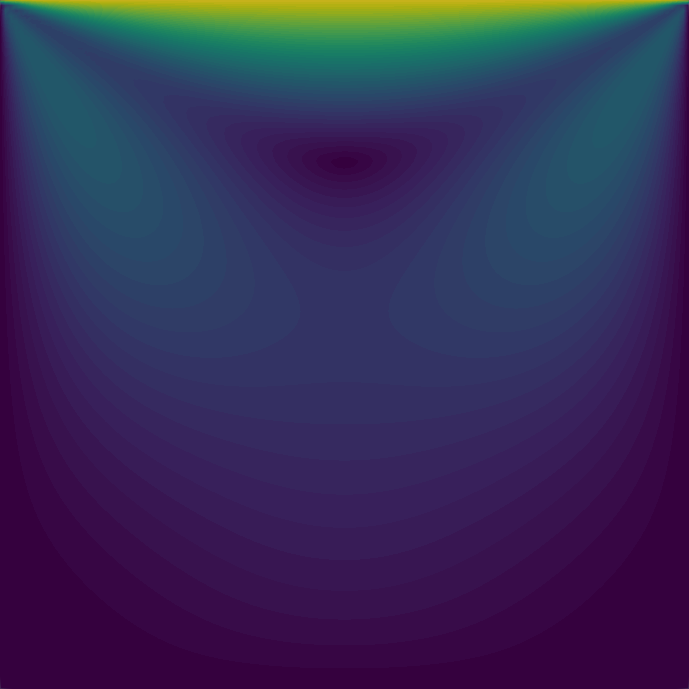}
	\caption{2d Driven Cavity: Geometry, coarse grid, and velocity field of the computed solution.}
	\label{img:cavity geometry}
\end{figure}

The computational domain is given by the unit square $\Omega = (0,1)^2$.
The flow is driven by a prescribed velocity 
{$u_D =  (1, 0)^T$}
on the top boundary $\Gamma_{top}$, 
and 
{$u_D =  (0, 0)^T$}
on the remaining part $\partial \Omega \backslash \Gamma_{top}$ 
of the boundary $\partial \Omega$ of $\Omega$; 
see Figure \ref{img:cavity geometry} for a visualization of this setting.
We first consider the stationary Stokes problem \eqref{eqn:Stokes:cf},
with the the viscosity $\mu = 1$ 
and the right-hand side $f = 0$.

\subsubsection{Different diagonal approximations of the Schur Complement}
\label{subsubsec:DrivenCavity2d:DifferentApproximations}

\begin{figure}[H]
	\centering
	\includegraphics[width=0.45\textwidth]{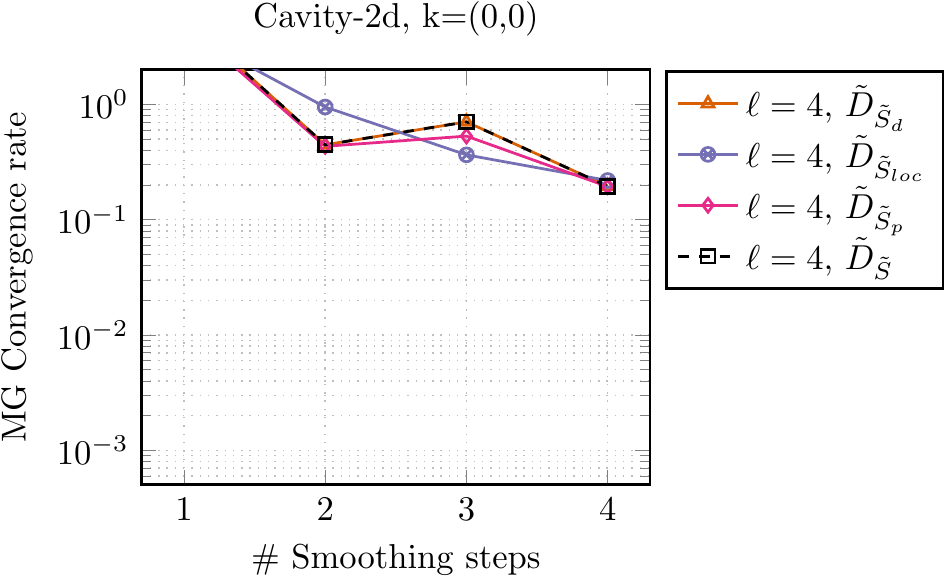}
	\includegraphics[width=0.45\textwidth]{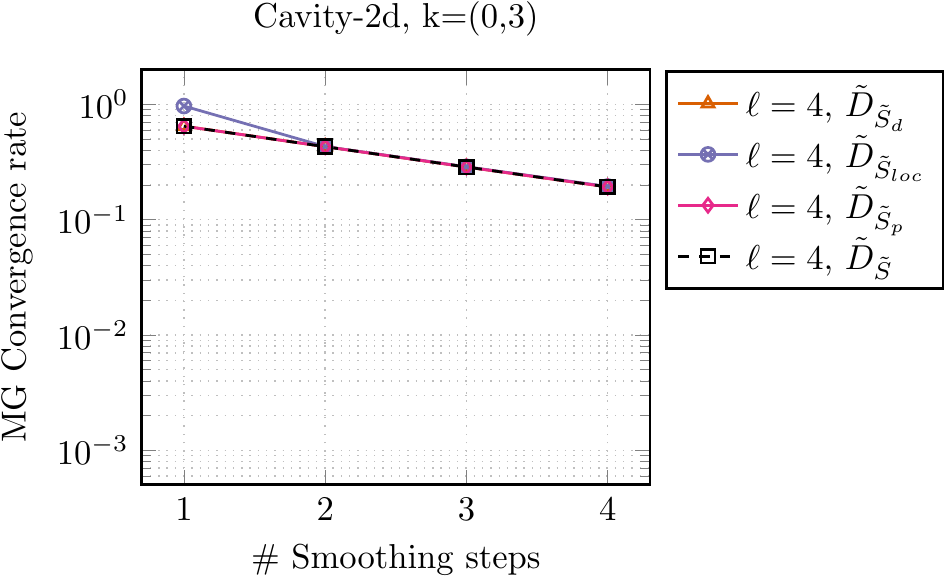}
	
	\vspace{1mm}
	
	\includegraphics[width=0.45\textwidth]{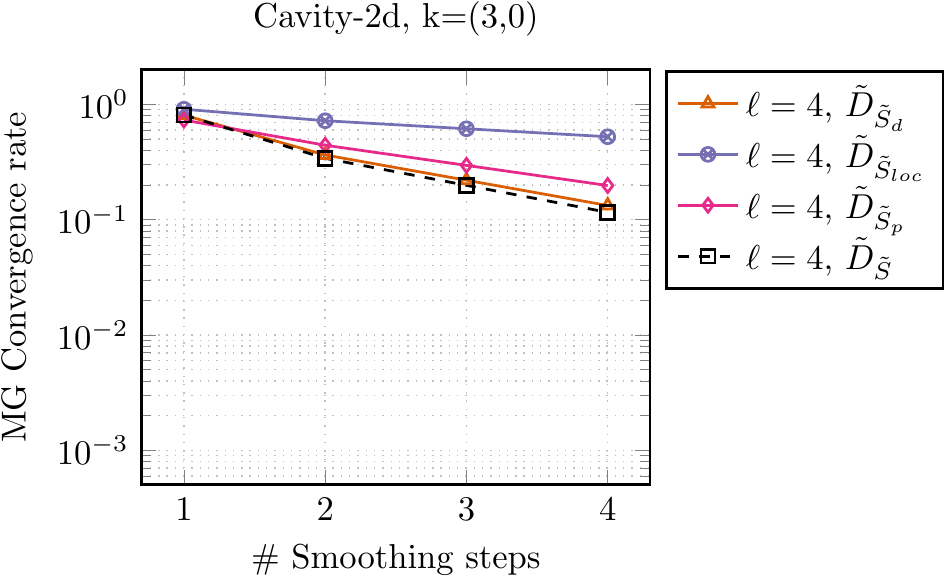}
	\includegraphics[width=0.45\textwidth]{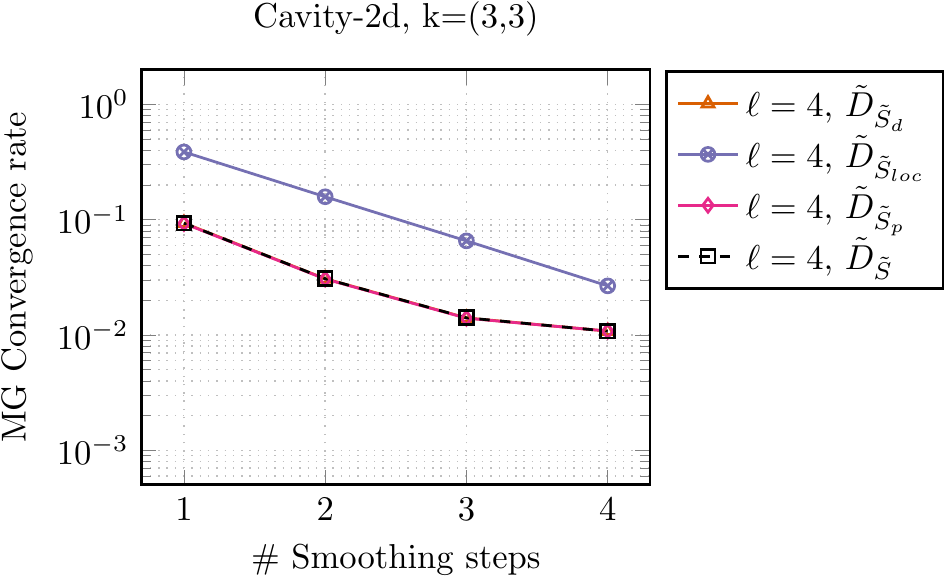}
	\caption{2d Driven Cavity: Convergence rates for different diagonal approximations of the Schur complement
	  and various Chebyshev polynomial degrees.}
	\label{fig:cavity:diagonal}
\end{figure}

Again, we investigate the effect of the approximations 
\eqref{eq:diag d}, \eqref{eq:diag p}, and \eqref{eq:diag loc} replacing the computationally expensive diagonal $D_{\tilde{S}}$.
	The convergence rates using the different approximations for $\tilde{D}_{\tilde{S}}$ are visualized in Figure \ref{fig:cavity:diagonal}.
We here consider only $2467$ dofs since we want to compare different cheaper approximations 
with the expensive diagonal $D_{\tilde{S}}$.
	The plot only shows the results for the Jacobi setting ($k=0$) and the Chebyshev-Jacobi degrees $k=3$, but the behavior is similar for other polynomial degrees in between,
	and will be investigated afterwards.
	We can see that using $\tilde{D}_{\tilde S_d}$ yields very similar rates as the exact diagonal $D_{\tilde S}$ for most degrees.
	Furthermore, using the pressure approximations works very well in this example.
	This might be due to the simple tensor-product geometry.
	In the Sch\"afer-Turek benchmark later on, this choice still works but is less effective.
	The locally computed diagonals fall slightly behind, most likely due to the missing contributions from neighboring elements.
	In any case, increasing the number of smoothing steps improves the convergence rate of the multigrid solver.
	In the Jacobi case 
	given in the top-left corner of Figure~\ref{fig:cavity:diagonal}, 
	this is not obvious from the plot.
	However, for sufficiently many smoothing steps, the expected behavior $q \sim 1/\sqrt{m}$ is attained.

\subsubsection{Different Polynomial Degrees}
\label{subsubsec:DrivenCavity2d:DifferentPolynomialDegrees}

\begin{figure}[htb]
	\centering
	\includegraphics[width=1.0\textwidth]{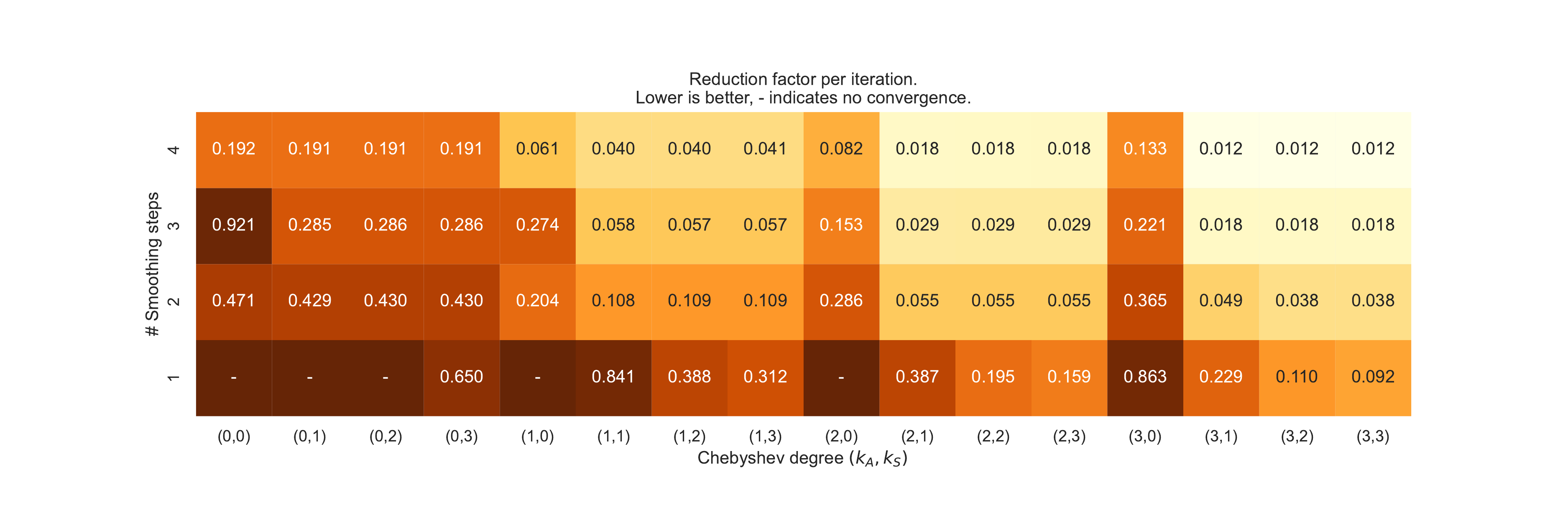}
	\caption{2d Driven Cavity: Convergence rates for varying Chebyshev degrees and smoothing steps.}
	\label{fig:cavity:rates}
\end{figure}

	Next, we look at the influence of the polynomial degrees $k_A$ and $k_S$ of the Chebyshev smoother.
	Higher values improve the reduction but come at additional costs in the application.
	The convergence rates are visualized in Figure \ref{fig:cavity:rates} for the diagonal approximation $\tilde{D}_{\tilde S_d}$.
	
	In Figure \ref{fig:cavity:rps}, 
	we additionally consider the time that 
	is needed for the application of 
	the smoother. 
	We note that
	increasing the degree of the Chebyshev polynomial may not improve the convergence rates enough to justify the increased computational effort.
	We present the achievable reduction of the residual per second for a fixed problem size of $148\, 739$ dofs.
	The convergence rates stay constant during $h$-refinement, but the cost of applying the smoother grows linearly with the number of dofs.
	A hyphen denotes configurations that did not yield a convergent algorithm, e.g., because an insufficient number of smoothing steps has been performed.
	We can see that the additional cost of higher Chebyshev degrees does not necessarily pay off.
	Rather, the best performance is achieved for degrees $(1,1), (2,1)$ with two smoothing steps and $(3,2)$ with one smoothing step.

\subsubsection{$W$-cycle versus $V$-cycle}	
\label{subsubsec:DrivenCavity2d:VversusW}

\begin{figure}[htb]
	\centering
	\includegraphics[width=1.0\textwidth]{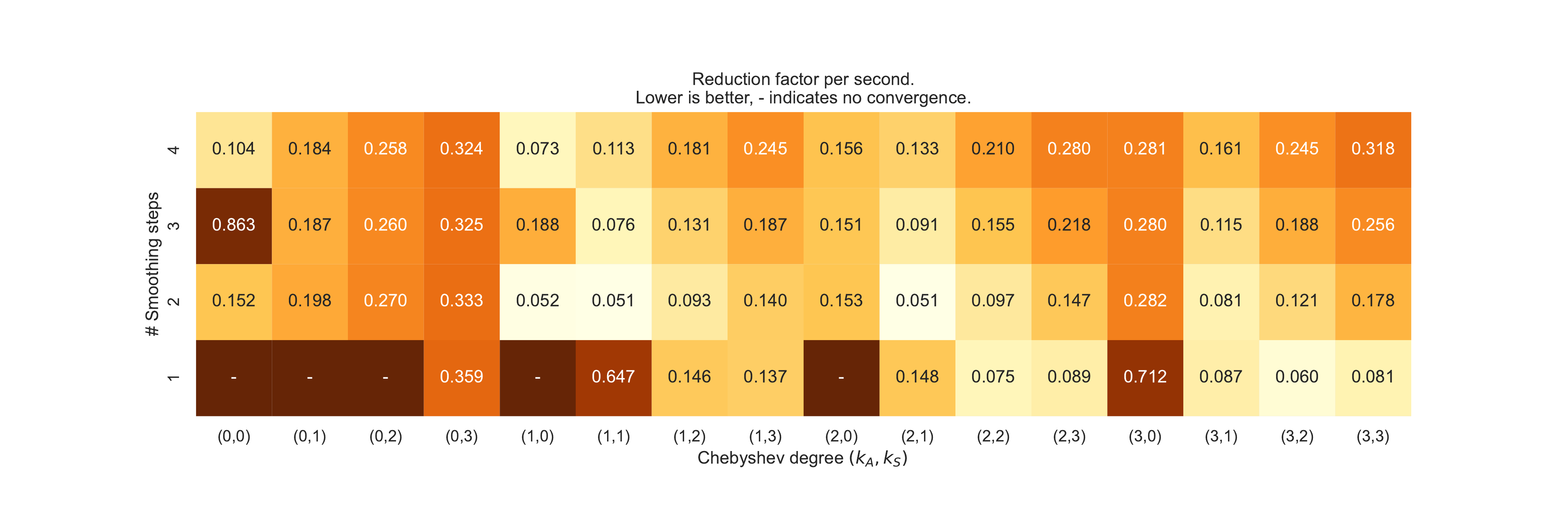}
	\includegraphics[width=1.0\textwidth]{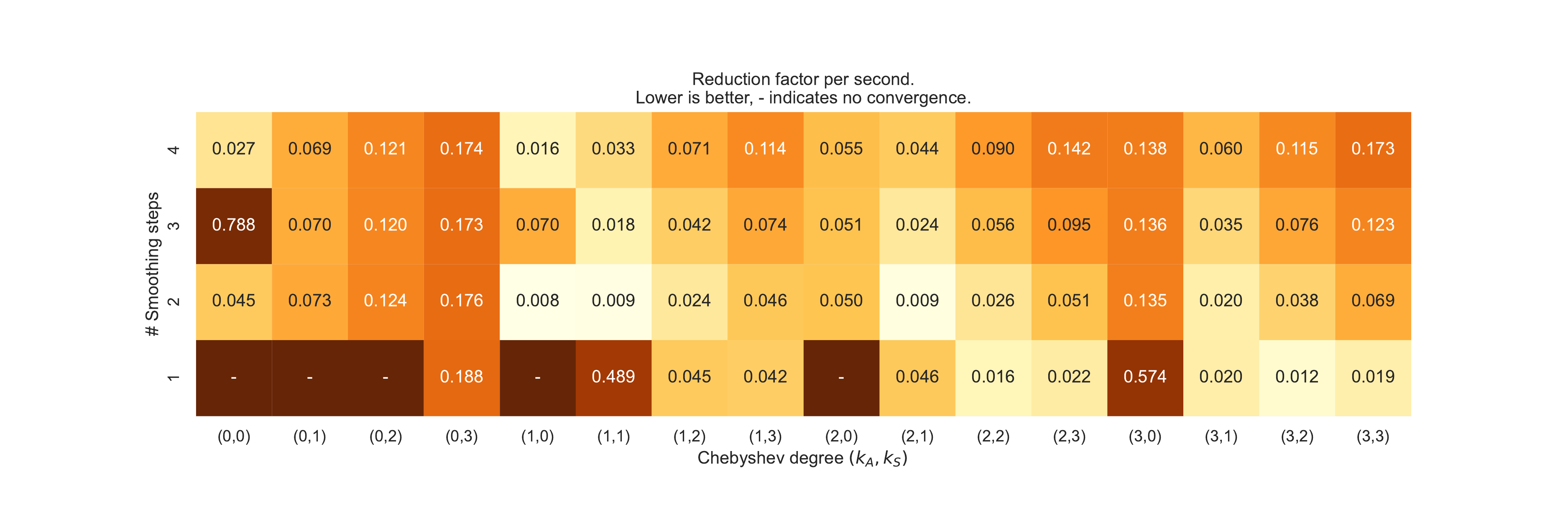}
	\caption{2d Driven Cavity: Reduction factor per second: (a) W-cycle, (b) V-cycle.}
	\label{fig:cavity:rps}
\end{figure}

The use of the $V$-cycle instead of the $W$-cycle yields almost the same convergence rates,
but the $V$-cycle is cheaper; see Figure \ref{fig:cavity:rps} (b).
To stay in line with the multigrid 
convergence analysis provided in Section~\ref{sec:MatrixfreeGMG},
we mainly use the $W$-cycle here.
However, for practical purposes, it is advisable to use a $V$-cycle for better performance.

\clearpage

\subsection{Sch\"afer-Turek Benchmark 2d-1}
\label{subsec:SchaeferTurekBenchmark}

In the third series of numerical tests, we apply the matrix-free geometrical multigrid solvers 
proposed and investigated in Section~\ref{sec:MatrixfreeGMG} 
to the  well-known Sch\"afer-Turek benchmark; see \cite{ScTuDu96}.

\subsubsection{Geometry, Boundary Conditions, and Parameters}
\label{subsubsec:SchaeferTurek:Data}

\begin{figure}[htb]
	\centering
	\includegraphics[width=0.7\textwidth]{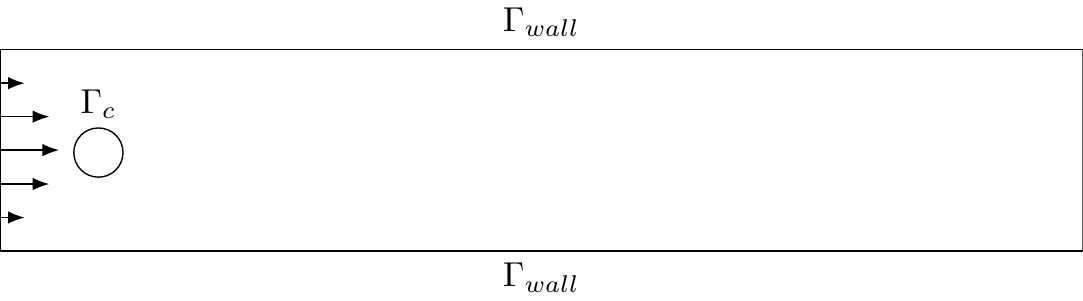}
	\vspace{1mm}
	\includegraphics[width=0.7\textwidth]{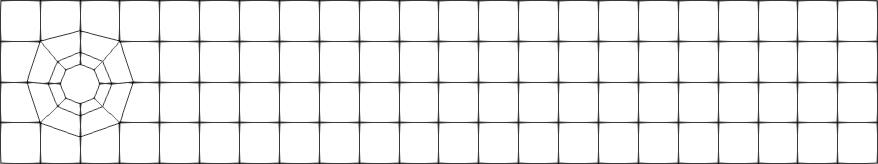}
	\vspace{1mm}
	\includegraphics[width=0.7\textwidth]{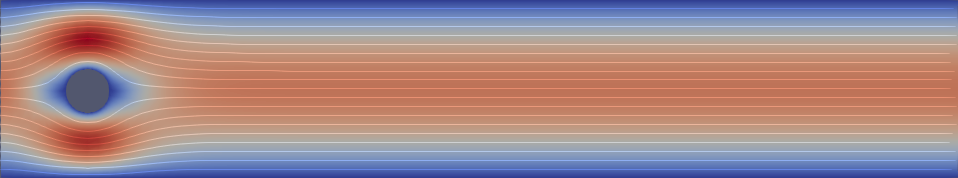}
	\vspace{1mm}
	\includegraphics[width=0.7\textwidth]{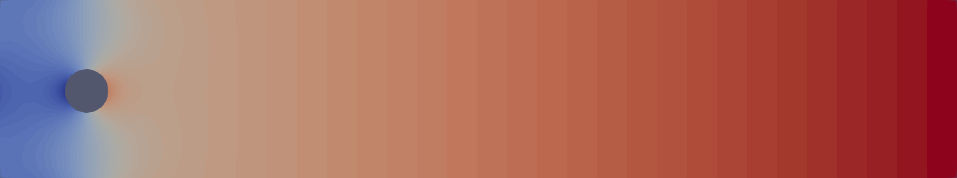}
	\caption{Sch\"afer-Turek: Geometry, coarse grid, velocity, and pressure field of the Sch\"afer-Turek benchmark.}
	\label{img:turek geometry}
\end{figure}

The flow inside the channel 
$(0, 2.2) \times (0, 0.41)$ around the circular obstacle that is located at $(0.2, 0.2)$ with diameter $0.1$
is driven 
by the given velocity
$u_D = (v_{in}(x_2), 0)$ 
on the left boundary $\Gamma_{in}$, 
with zero volume forces $f = 0$.
The inflow profile 
is prescribed as
$v_{in}(x_2) = 4\, v_\text{max}\, x_2 \,(0.41 - x_2) (0.41)^{-2}$ 
with the maximum velocity $v_{max} = 1$.
The velocity is zero on the top and bottom parts $\Gamma_{wall}$, as well as on the cylinder boundary $\Gamma_{c}$ (no-slip conditions).
Figure~\ref{img:turek geometry} provides an illustration of the computational geometry, the mesh,
and the computed velocity and the pressure fields.
This example is again stationary ($\gamma = 0$) with viscosity $\mu = 1$.

\subsubsection{Different Diagonal Approximations of the Schur Complement}
\label{subsubsec:SchaeferTurek:DifferentApproximations}

\begin{figure}[htb]
	\centering
	\includegraphics[width=0.45\textwidth]{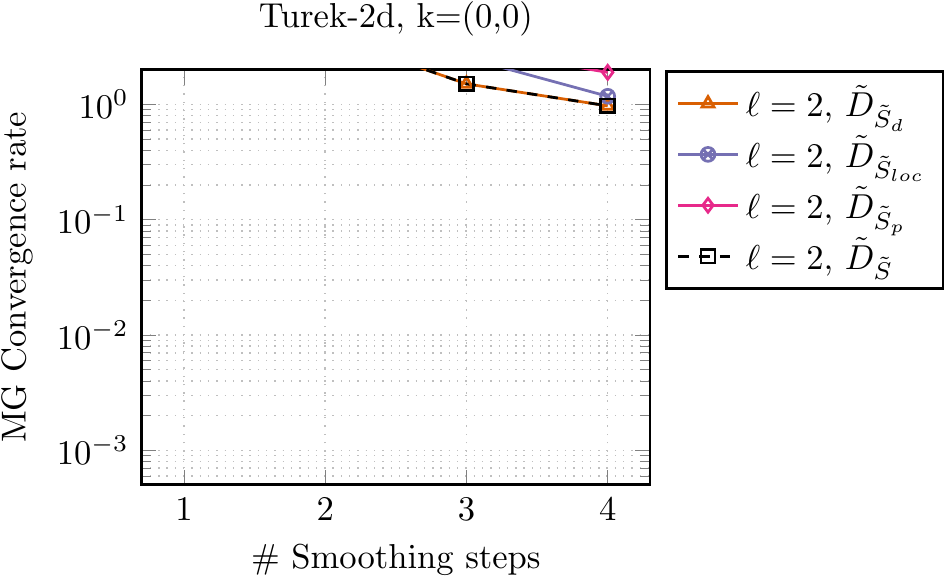}
	\includegraphics[width=0.45\textwidth]{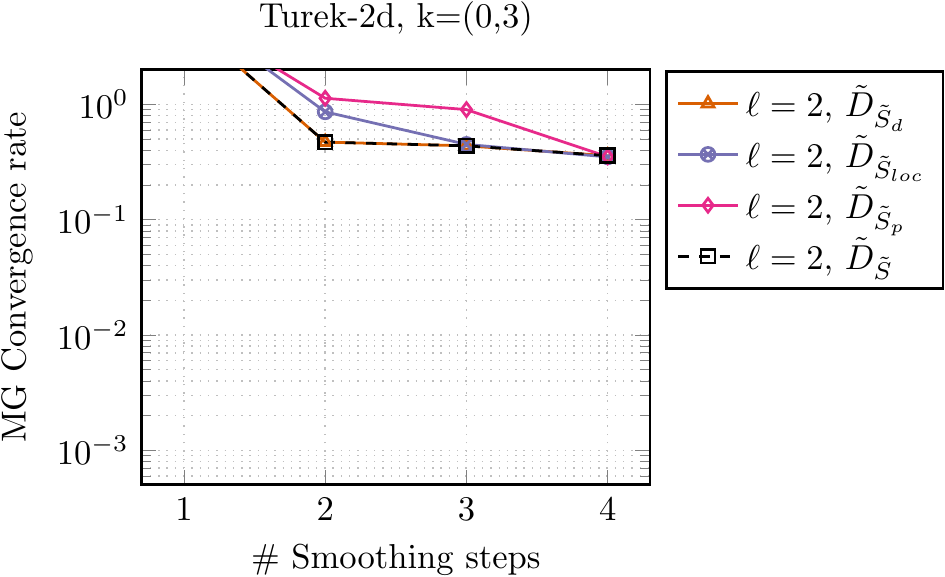}
	
	\vspace{1mm}
	
	\includegraphics[width=0.45\textwidth]{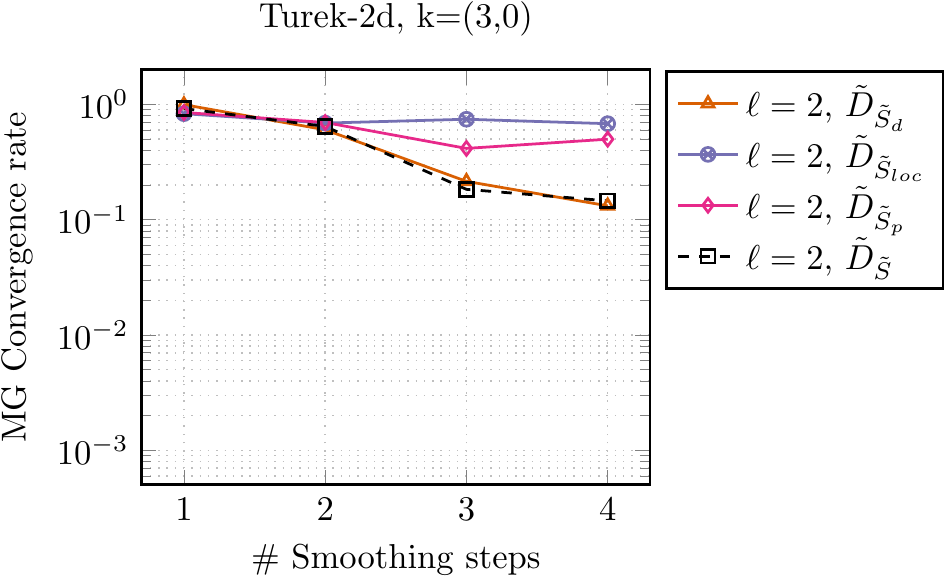}
	\includegraphics[width=0.45\textwidth]{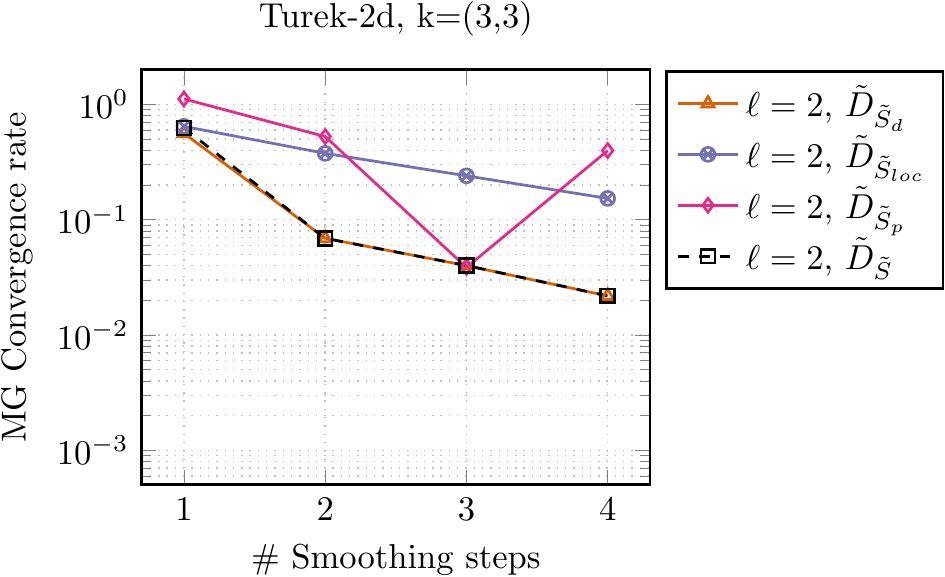}
	\caption{Sch\"afer-Turek: Convergence rates for different diagonal approximations and Chebyshev degrees.}
	\label{fig:turek:diagonal}
\end{figure}

As in the first example, we first consider a 
relatively coarse discretization
with $3900$ fine grid dofs, and compare the various diagonal approximations \eqref{eq:diag d}, \eqref{eq:diag p}, and \eqref{eq:diag loc},  in Figure \ref{fig:turek:diagonal}.
We can see that four smoothing steps are not enough to yield a convergent multigrid method in the Jacobi case $k_A = k_S = 0$.
Similar to the example before, the diagonal approximation $\tilde{D}_{\tilde S_d}$ yields almost the same rates as the exact variant $D_{\tilde S}$.
The locally assembled diagonal $\tilde{D}_{\tilde S_{loc}}$ results in consistently worse convergence rates.
In this example, the pressure approximation does not work as well as in the 2d driven cavity test case, with the surprising exception of $k=(3,3)$ and three smoothing steps.
So far, we have no explanation for this behavior.

\subsubsection{Different Polynomial Degrees and more Unknowns}
\label{subsubsec:SchaeferTurek:DifferentPolynomialDegrees}

\begin{figure}[H]
	\centering
	\includegraphics[width=1.0\textwidth]{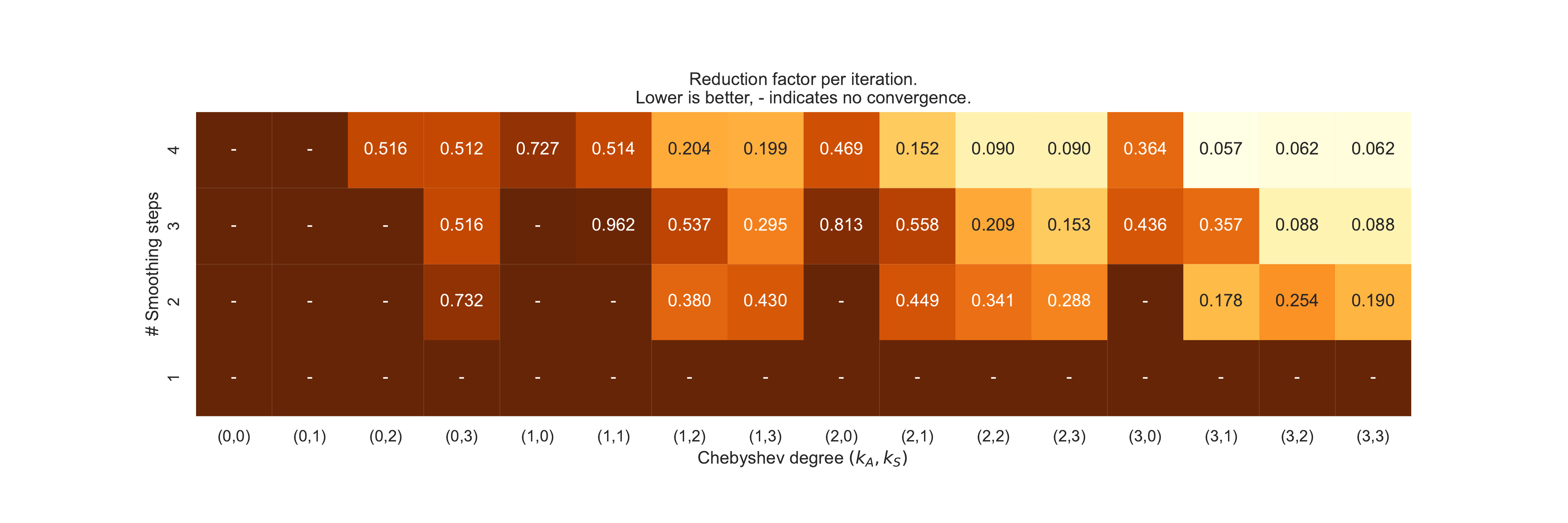}
	\caption{Sch\"afer-Turek: Convergence rates for varying Chebyshev degrees and smoothing steps.}
	\label{fig:turek:rates}
\end{figure}

\begin{figure}[H]
	\centering
	\includegraphics[width=1.0\textwidth]{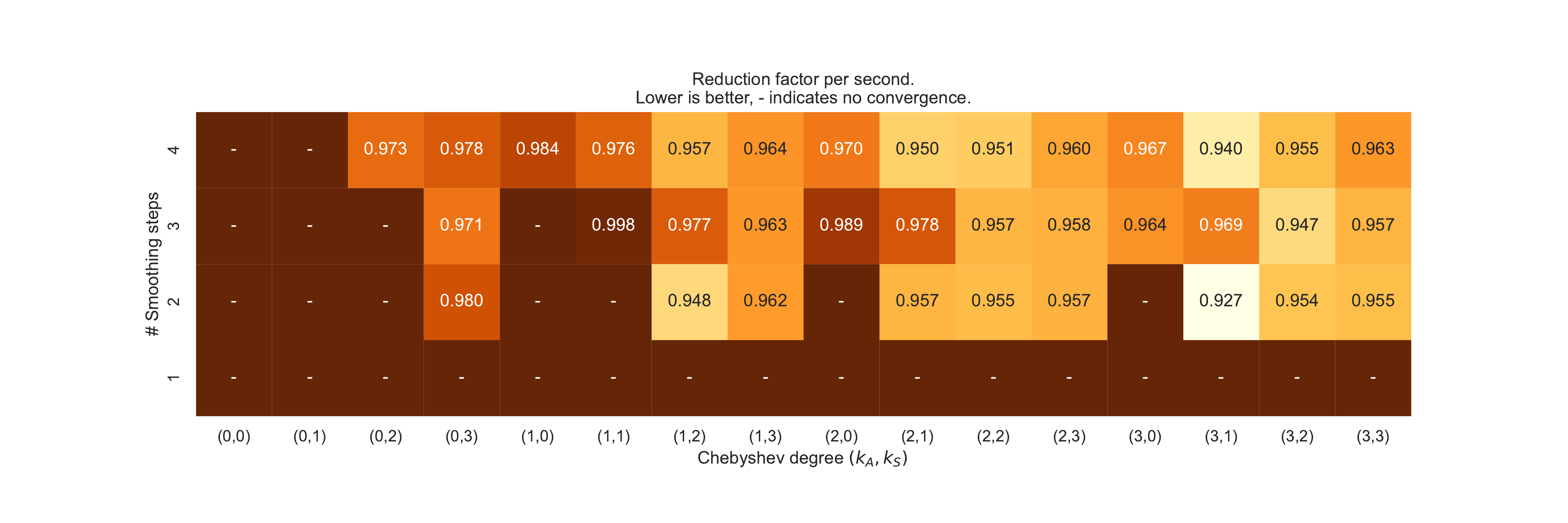}
	\caption{Sch\"afer-Turek: Reduction factor per second.}
	\label{fig:turek:rps}
\end{figure}

In the case of more refinement levels resulting in $232\, 800$ fine grid dofs,
the convergence rates  for all tested combinations of $k_A, k_S,$ and $m$ are 
visualized in the heatmap shown in Figure \ref{fig:turek:rates}.
In contrary to the 2d driven cavity example, 
a single smoothing step is not enough to guarantee convergence even in the case $k_A = k_S = 3$.
From Figure \ref{fig:turek:rps}, we deduce that two smoothing steps with degrees $k_A=3$ and $k_S=1$ yield the 
best reduction per second, 
and hence, the fastest overall solver.

\subsection{3d Driven Cavity}
\label{subsec:DrivenCavity3d}
In the fourth series of numerical tests,
we illustrate the performance of the smoother on a $3d$ extension of the lid-driven cavity benchmark 
from  Subsection~\ref{subsec:DrivenCavity2d}.
\subsubsection{Geometry, Boundary conditions, and Parameters}
\label{subsubsec:DrivenCavity3d:Data}

\begin{figure}[H]
	\centering
	\includegraphics[width=0.2\textwidth]{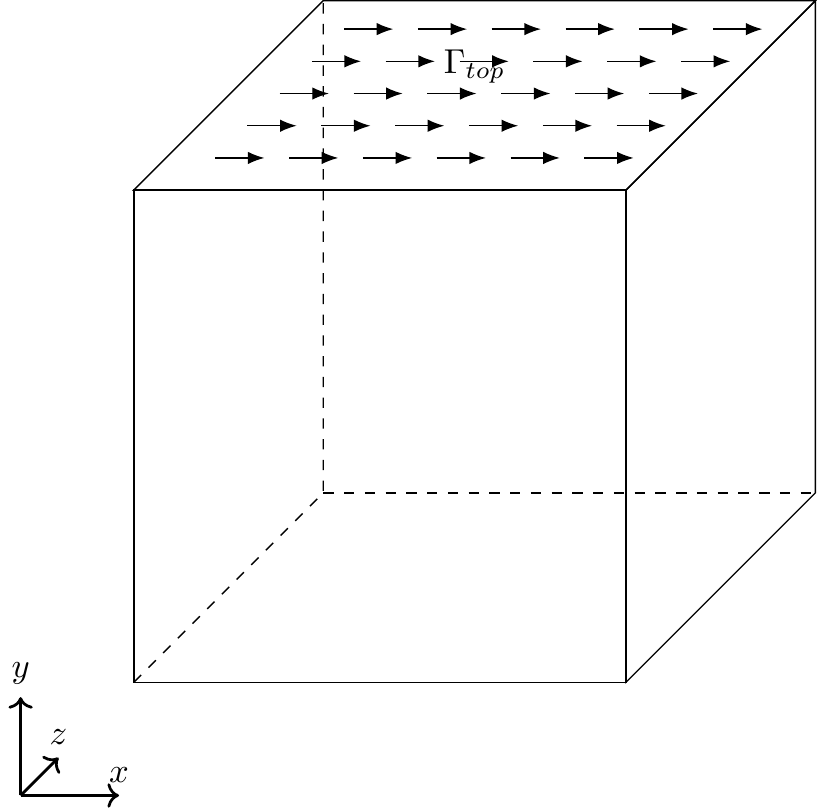}
	\hspace{1cm}
	\includegraphics[width=0.2\textwidth]{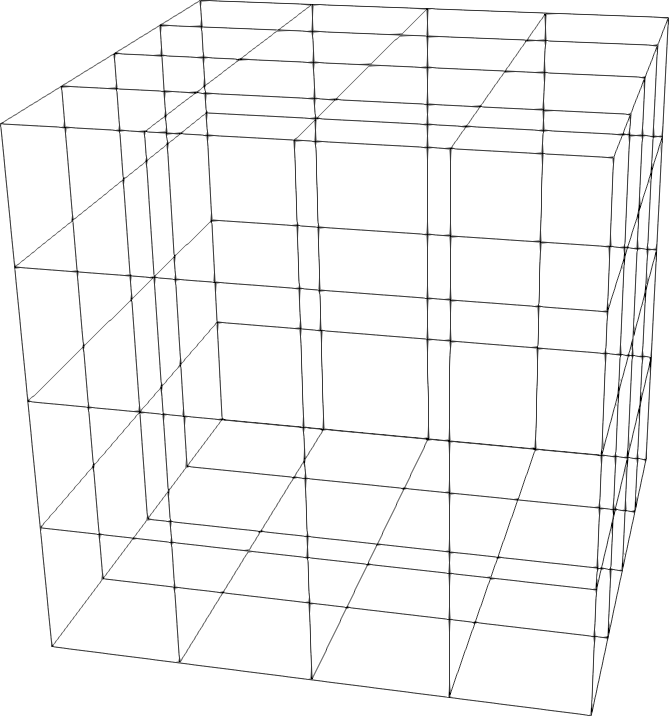}
	\hspace{1cm}
	\includegraphics[width=0.2\textwidth]{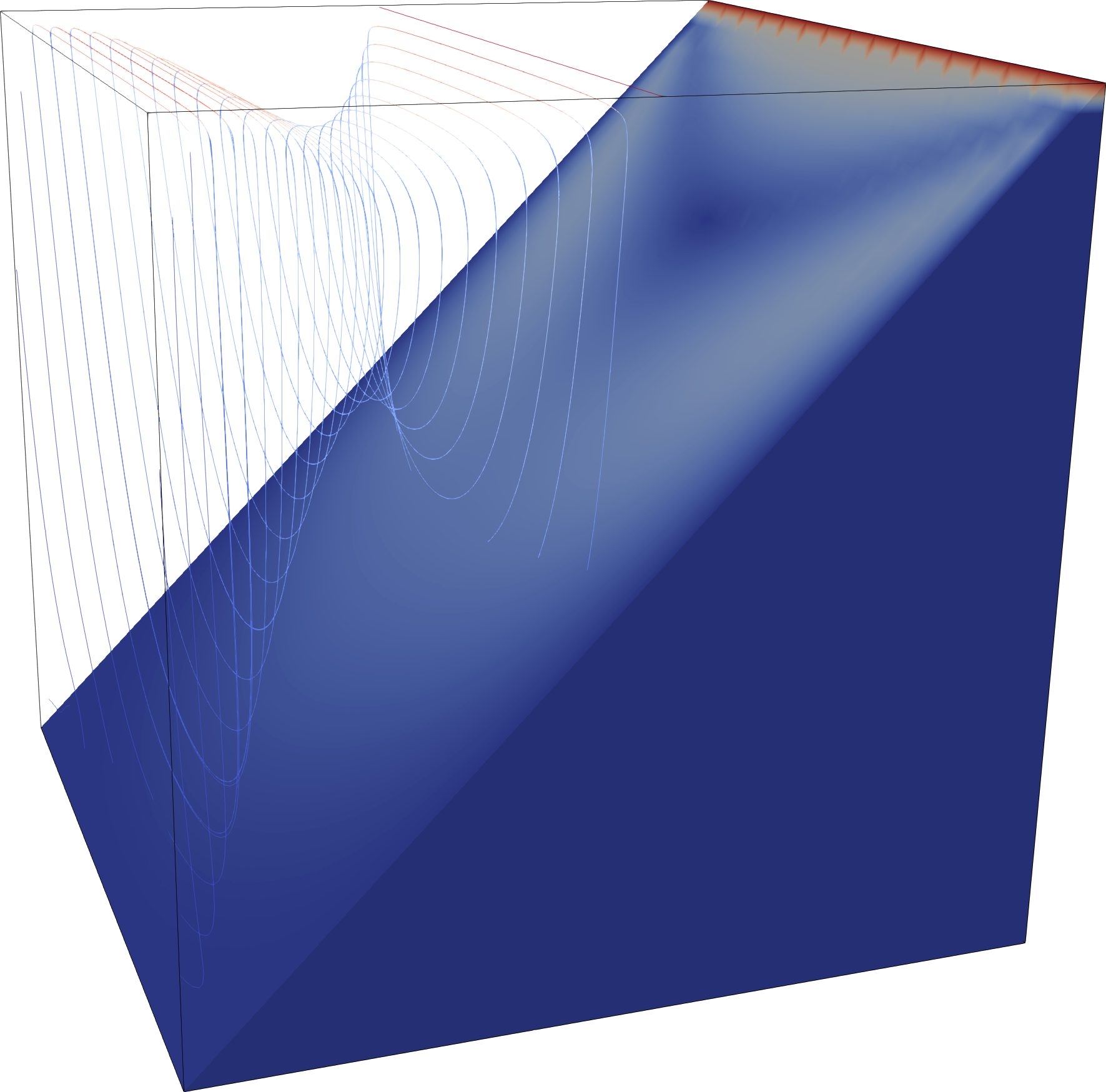}
	
	\caption{3d Driven Cavity: Geometry, coarse grid, and velocity field of the solution for the Lid-Driven Cavity benchmark in $3d$.}
	\label{img:cavity3d:geometry}
\end{figure}

The computational domain is 
nothing but the unit cube
$\Omega = (0,1)^3$.
The fluid is driven by the constant velocity $u_D = (1,0,0)$ 
on the top boundary $\Gamma_{top} = (0,1) \times (0,1) \times \{ 1 \}$.
On the remaining boundary we prescribe $u_D = 0$.
Figure \ref{img:cavity3d:geometry} illustrates the geometry, the mesh, and 
the computed solution 
for the stationary Stokes problem \eqref{eqn:Stokes:cf}, with $\mu = 1$ and $f = 0$.

\begin{figure}[H]
	\centering
	\includegraphics[width=0.45\textwidth]{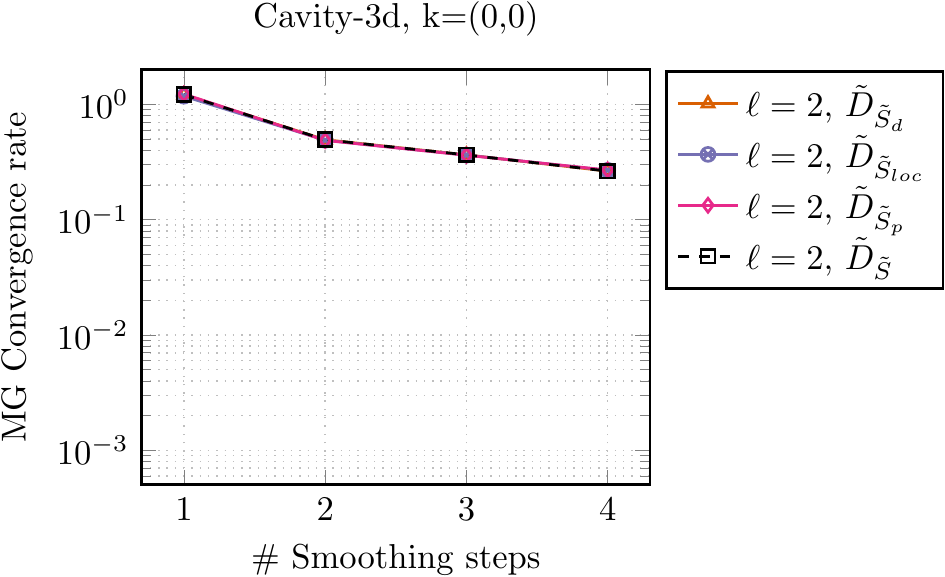}
	\includegraphics[width=0.45\textwidth]{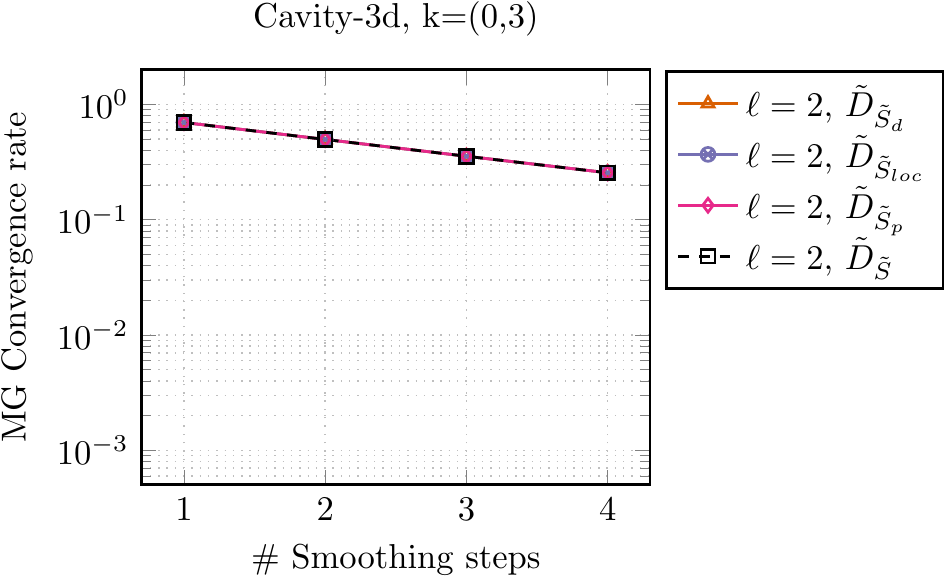}
	
	\vspace{1mm}
	
	\includegraphics[width=0.45\textwidth]{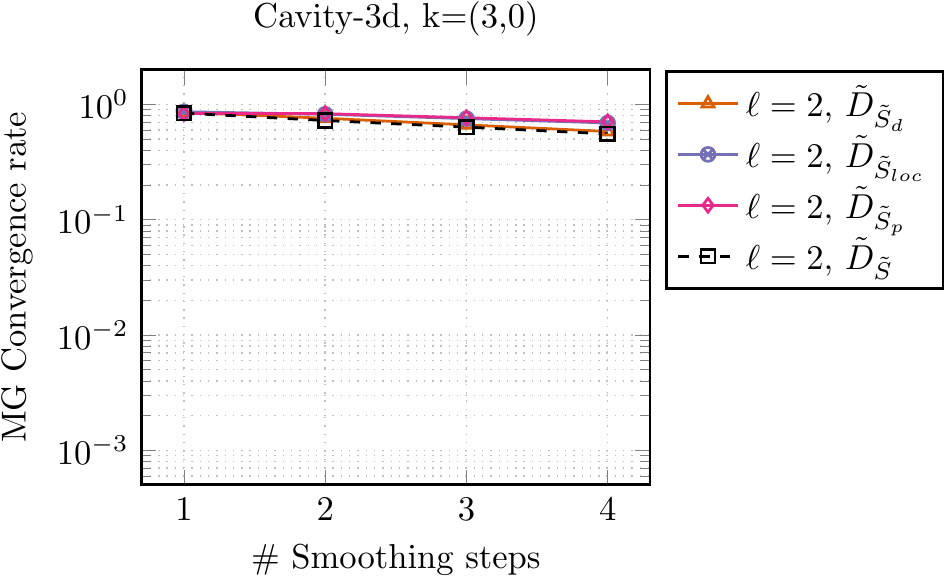}
	\includegraphics[width=0.45\textwidth]{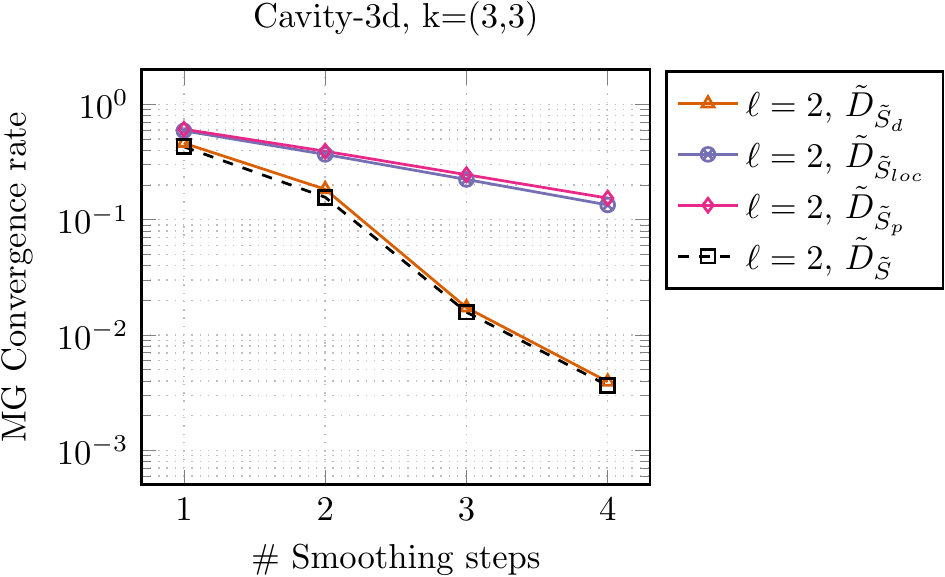}
	
	\caption{3d Driven Cavity: Convergence rates for different diagonal approximations and various Chebyshev degrees.}
	\label{fig:cavity3d:diagonal}
\end{figure}

\subsubsection{Discussion of the Numerical Behavior of the Solver}
\label{subsubsec:DrivenCavity3d:Solver}

\begin{figure}[H]
	\centering
	\includegraphics[width=1.0\textwidth]{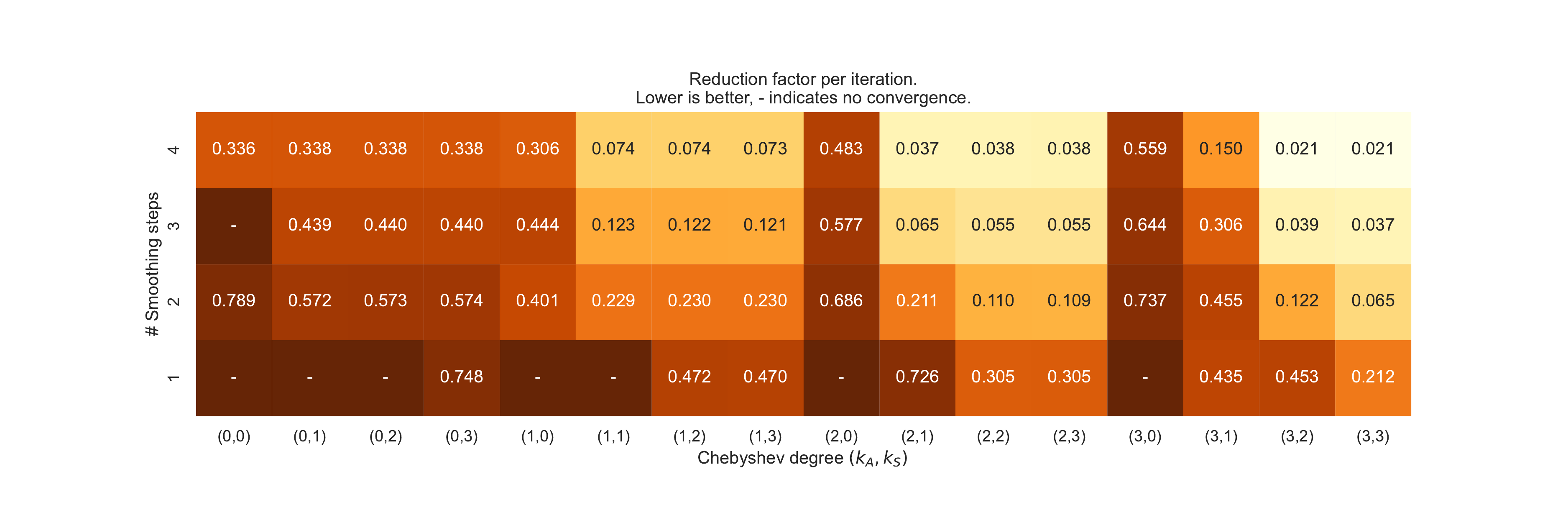}
	\caption{3d Driven Cavity: Convergence rates for varying Chebyshev degrees and smoothing steps.}
	\label{fig:cavity3d:rates}
\end{figure}

\begin{figure}[H]
	\centering
	\includegraphics[width=1.0\textwidth]{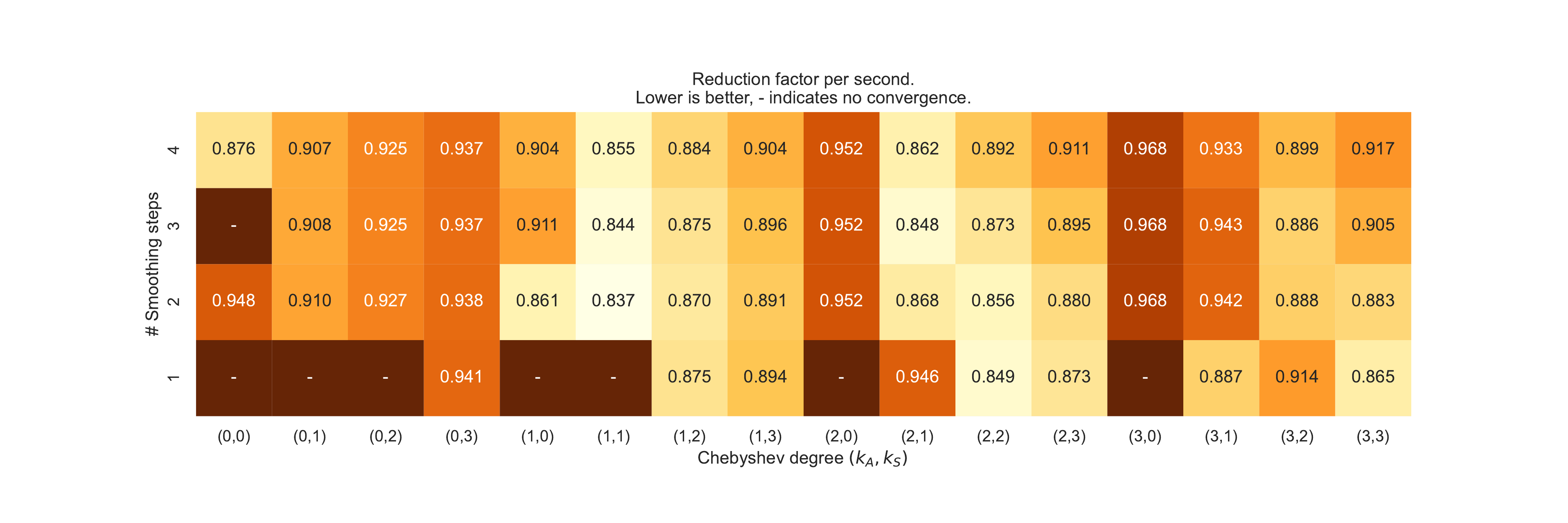}
	\caption{3d Driven Cavity: Reduction factor per second.}
	\label{fig:cavity3d:rps}
\end{figure}

For the computations, we used a problem size of $2312$ dofs for the comparison of the diagonals, and $112\, 724$ dofs for the convergence rate and timing test.
	The results in $3d$ are in line with previous observations
as shown in the Figures \ref{fig:cavity3d:rates} and \ref{fig:cavity3d:rps}.
	Interestingly, $k=(3,0)$ yields worse rates compared to $k=(0,3)$ or the plain Jacobi setting, which was not the case in the 2d example.
	The highest reduction per second is again achieved with two smoothing steps of the $(1,1)$ Chebyshev smoother.

\subsection{Generalized Stokes Problem}
\label{subsec:GeneralizedStokesProblem}
The last series of numerical experiments is devoted to the numerical solution 
of the generalized Stokes problem~\eqref{eqn:genStokes:cf} arising from 
the implicit Euler discretization of the instationary Stokes equations
for the $2d$ driven cavity problem.

\subsubsection{Configuration}
\label{subsubsec:GeneralizedStokesProblem:Configuration}
The basic setting is the same as in Subsection~\ref{subsec:DrivenCavity2d}.
In addition, we now work with $\gamma = 1/\Delta t$ 
arising from the time discretization.
Moreover, we consider time-step sizes $\Delta t = 1$ and $\Delta t = 0.01$, 
with homogeneous initial conditions at $t=0$.
We only present the results for a single time step here. However, the convergence rates stay the same for longer simulations as expected.

\subsubsection{Discussion of the Numerical Behavior of the Solver}
\label{subsubsec:GeneralizedStokesProblem:Solver}

\begin{figure}[H]
	\centering
	\includegraphics[width=1.0\textwidth]{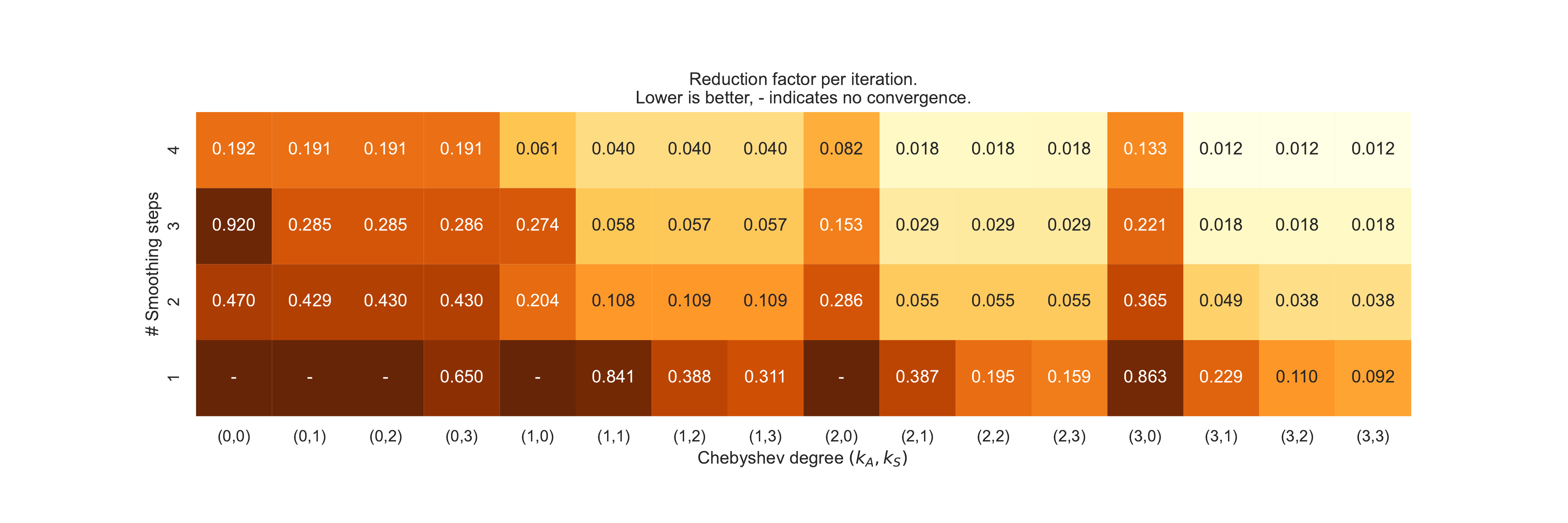}
	\includegraphics[width=1.0\textwidth]{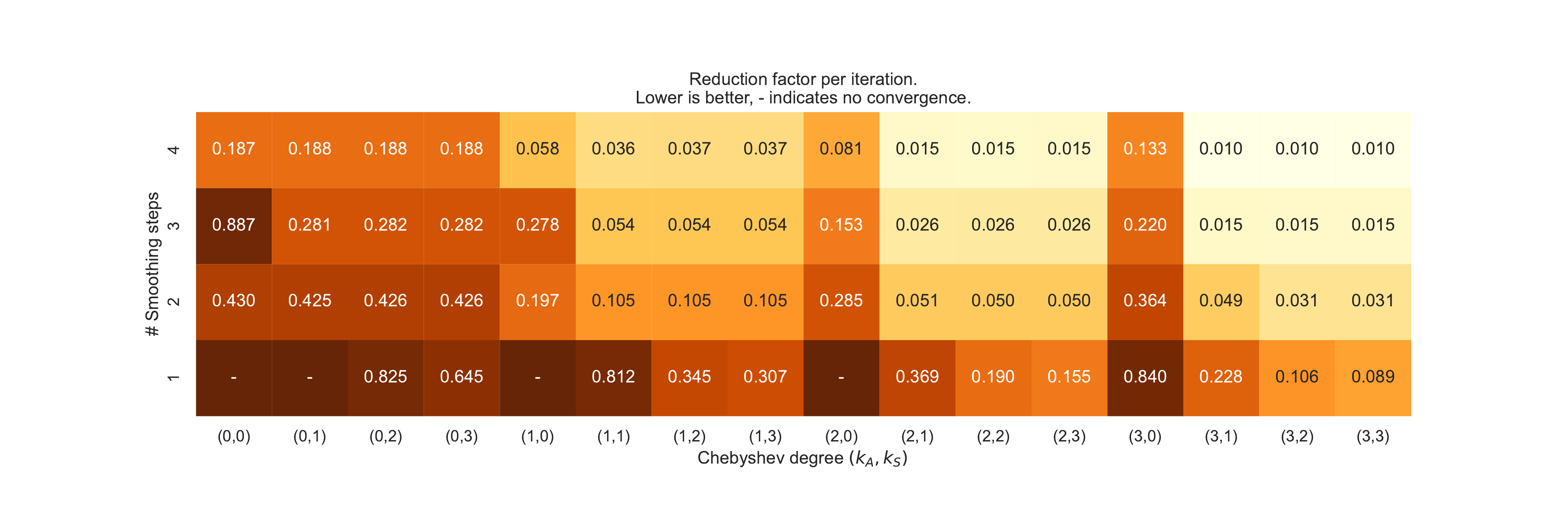}
	\caption{Time-dependent 2d Driven Cavity Problem: Convergence rates for varying Chebyshev degrees and smoothing steps. (a) $\Delta t = 1$, (b) $\Delta t = 0.01$.}
	\label{fig:dt:rates}
\end{figure}

From Figure \ref{fig:dt:rates}, we can see that the convergence rates for $\Delta t = 1$ and $\Delta t = 0.01$ are almost indistinguishable from each other, and also match with the results in Figure \ref{fig:cavity:rates} ($\Delta t = \infty$).
Thus, also the reduction per second is the same, 
since the additional mass matrix stemming from the time discretization does not impose a significant overhead.

\section{Conclusions and Outlook}
\label{sec:ConclusionsOuilook}
We proposed and analyzed matrix-free, monolithic geometric multigrid solvers that work 
for Stokes and generalized Stokes problems as the numerical studies performed for different
well-known benchmark problems showed. In particular, the Chebyshev acceleration 
in $\hat{A}$ and $\hat{S}$ improves
the smoothing property of the symmetric inexact Uzawa smoother \eqref{eqn:isus1} - \eqref{eqn:isus3}
considerable.  
Indeed, increasing the polynomial degree of the Chebyshev smoother improves the convergence rates significantly, with moderate increase in the computational effort.
	The best performance requires a fine tuning between the number of smoothing steps and the degrees of the Chebyshev smoothers.
	However, choosing a high degree (e.g., $3$) yields very similar performance and is a good default initial choice.
Although we did not present parallel studies in this paper, the code runs efficiently in parallel with excellent strong scalability.
The parallelization is done without changes in the results beyond floating point errors.
Hence, all convergence rates are independent of the number of utilized cores.

The Stokes and generalized Stokes solvers discussed in this paper can also be used for 
solving stationary and instationary Navier-Stokes problems, and as building block in 
solvers or preconditioners for 
other coupled problems involving flow, such as fluid-structure interaction or Boussinesq approximations 
for laser-beam 
material processing or convection of material in the earth's mantle.
However, then we have to treat the 
convection term on the right-hand side. 
So, we cannot expect robustness of the solver with respect to dominating convection.
The Oseen or the Newton linearizations that takes care of the convection 
also lead to systems of the form \eqref{eqn:LinearSystemKx=f},
but with an non-symmetric and, in general, non-positive definite matrix $A$.
This changes the game completely. The theory for \cite{SchoeberlZulehner:2003a} is not applicable anymore.
It is certainly a challenge to develop a theory that allows us to analyze this case.
One possible approach is to use Reusken's lemma, which allows to study the smoothing property also in the non-symmetric case, see \cite{Reusken_1992,EckerZulehner_1996}.

\section*{Appendix}
\appendix

\section{Discrete inf-sup Conditions}

Here we assume that $\mathcal{T}_h$ is a subdivison of $\Omega$ such that each $\tau \in \mathcal{T}_h$ satisfies the following conditions.
\begin{description}
\item[In 2d:]
$\tau = T_\tau((0,1)^2)$ is a convex quadrilateral, $T_\tau$ bilinear, and $\overline{\tau}$ has at most two edges on $\partial \Omega$, among them there is no pair of opposite edges.
\item[In 3d:] $\tau = T_\tau((0,1)^3)$ is a parallelepiped, $T_\tau$ affine linear, and $\overline{\tau}$ has at most three faces on $\partial \Omega$, among them there is no pair of opposite faces.  
\end{description}
Furthermore, we assume that the standard compatibity and shape regularity conditions hold for $\mathcal{T}_h$.

Then we have for the Taylor-Hood $\mathcal{Q}_k$-$\mathcal{Q}_{k-1}$ finite element the following results.
\begin{Theorem}
\label{discreteInf-Sup}
Let $k \in \mathbb{N}$ with $k\ge 2$.
Then there is a constant $c > 0$ such that
\[
  \sup_{0 \neq v_h \in V_h}
  \frac{\displaystyle b(v_h,q_h)}{\|v_h\|_{H^1(\Omega)^d}}
    \ge c \, \|q_h\|_{L^2(\Omega)} 
        \quad \text{for all} \ q_h \in Q_h .
\]
\end{Theorem}
This is equivalent to the classical inf-sup condition \eqref{eqn:infsup}.

The next inf-sup condition has its origin in \cite{BercovierPironneau_1979} and is formulated for a different pair of norms. 

\begin{Theorem} \label{BerPirglobal}
Let $k \in \mathbb{N}$ with $k\ge 2$.
Then there is a constant $c > 0$ such that for all $\tau \in \mathcal{T}_h$
\[
  \sup_{0 \neq v_h \in V_h}
  \frac{\displaystyle b(v_h, q_h)}{\|v_h\|_{L^2(\Omega)^d}}
    \ge c \, \|\nabla q_h\|_{L^2(\Omega)^d} 
        \quad \text{for all} \ q_h \in Q_h.
\]
\end{Theorem}

The proofs in \cite{BercovierPironneau_1979} were restricted to meshes of rectangles in 2d resp.~bricks in 3d.

We need also the following local version of Theorem \ref{BerPirglobal}.

\begin{Theorem} 
\label{BerPirlocal}
Let $k \in \mathbb{N}$ with $k\ge 2$.
Then there is a constant $c > 0$ such that
\[
  \sup_{0 \neq v_h \in V_\tau}
  \frac{\displaystyle b_\tau(v_h,q_h)}{\|v_h\|_{L^2(\tau)^d}} \ge c \, \|\nabla q_h\|_{L^2(\tau)^d}
        \quad \text{for all} \ q_h \in Q_\tau 
\]
with
\[
    b_\tau(v,q) = -\int_\tau q \, \div v \ d x, \quad
    V_\tau = \{ v|_\tau \colon v \in V_h \}, \quad
    Q_\tau = \{ q|_\tau \colon q \in Q_h \}.
\]
\end{Theorem}

The last two theorems translate to the following estimates in matrix notation.

\begin{equation}
\label{BerPirlocalMatrixVersion}
  B M_{v}^{-1} B^\top \ge c^2 \,  K_{p}
  \quad \text{and} \quad
  B_\tau M_{v,\tau}^{-1} B_\tau^\top \ge c^2 \, K_{p,\tau}.
\end{equation}

For proofs of Theorems \ref{discreteInf-Sup}, \ref{BerPirglobal}, and \ref{BerPirlocal}, 
see \cite{Zulehner:2022a}
For a proof of Theorem \ref{discreteInf-Sup} in 2d see also \cite{Stenberg1990}.

\bibliographystyle{abbrv}
\bibliography{JodlbauerLangerWickZulehner}
\end{document}